\documentclass[amsppt,11pt]{amsart}

\usepackage{amsmath,amsthm,amssymb}
\usepackage{verbatim}
\allowdisplaybreaks[4]
\title[Scaling limits of UST]
{Scaling limits of the uniform spanning tree and loop-erased random 
walk on finite graphs}
\author[Yuval Peres\,\, David Revelle]
{Yuval Peres$^*$, David Revelle$^\dagger$}
\date{June 6, 2005.
\newline\indent
$^*$Research supported by NSF grants
\#DMS-0104073 and \#DMS-0244479.
\newline\indent
$^\dagger$Research partially supported by an NSF postdoctoral
fellowship}
\def\eps{\varepsilon}
\def\R{\mathbb {R}}
\def\Z{\mathbb {Z}}

\def\F{\mathbb {F}}
\def\E{\mathbb{E}}

\def\P{\mathbb{P}}

\def\limn{\lim_{n\rightarrow \infty}}

\def\1{{\mathbf 1}}
\def\0{{\mathbf 0}}

\newtheorem{theorem}{Theorem}[section]
\newtheorem{lemma}{Lemma}[section]

\newtheorem{definition}{Definition}
\newtheorem{corollary}{Corollary}[section]

\def\concat{\ast}
\def\Cap{\mathop{\mathrm{Cap}_r}}
\newcommand{\Capr}[1]{\mathop{\mathrm{Cap}_{#1}}}

\def\LE{\mathbf{LE}}
\def\LLE{\mathbf{LE}_s}
\def\ind{\mathrm{I}}
\def\Iij{\ind_{ij}}

\def\tildeIij{\tilde \ind_{ij}}

\def\Ikj{\ind_{kj}}
\def\Ikl{\ind_{k\ell}}
\def\tildeIkl{\tilde \ind_{k \ell}}
\def\tildeIkj{\tilde \ind_{k j}}
\def\vol{|G_n|}
\def\sqrtvol{\vol^{1/2}}
\def\delsqrtvol{\vol^{1/2-\delta}}
\def\T{\mathcal{T}}
\def\Told{T_{\mathrm{old}}}
\def\Tnew{T_{\mathrm{new}}}
\def\killT{\widehat T}
\def\hatY{\widehat Y}
\def\hatZ{\widehat Z}
\def\hatW{\widehat W}
\def\cc{\theta}
\newcommand{\step}[1]{\mathcal{A}_{#1}}
\newcommand{\good}[1]{\mathcal{G}_{#1}}
\newcommand{\single}[1]{\mathcal{C}_{#1}}
\newcommand{\bad}[1]{\mathcal{B}_{#1}}
\newcommand{\Kindex}[1]{\widetilde{\mathcal{S}}_{#1}}
\newcommand{\Gindex}[1]{\mathcal{S}_{#1}}
\newcommand{\close}[2]{\mathop{\mathrm{Close}(#1, #2)}}
\newcommand{\lclose}[2]{\close{\LLE(A_{#1})}{ \LLE(A_{#2})}}
\newcommand{\extG}[2]{G_{#1,#2}^*}
\newcommand{\extK}[2]{K_{#1,#2}^*}
\def\mix{\tau}
\def\mixn{\mix_n}   
\def\F{\widehat \T}
\def\extT{\T^*}
\newcommand{\dT}[2]{d_\T(#1,#2)}
\newcommand{\dF}[2]{d_{\F}(#1, #2)}

\newcommand{\dTk}[2]{d_{\T_k}(#1, #2)}
\newcommand{\dTktwo}[2]{d_{\T_{k+1}}(#1,#2)}
\newcommand{\dtildeTk}[2]{d_{\widetilde\T_k}(#1, #2)}

\begin{document}

\begin{abstract}
Let $x$ and $y$ be chosen uniformly in a graph $G$.
We find the limiting distribution of the length of a 
loop-erased random walk from $x$ to $y$ on a 
large class of graphs that include the torus
$\Z_n^d$ for $d\geq 5$.  Moreover, on this family
of graphs we show that a suitably normalized finite-dimensional
scaling limit of the uniform spanning tree is
a Brownian continuum random tree. 
\end{abstract}

\maketitle

\section{Introduction}
\label{introduction}

The {\it uniform spanning tree} (UST) $\T$ on a graph $G$ is a 
random tree, uniformly distributed among all 
spanning trees of $G$.  
For two points $x$ and $y$ in a graph, let $\dT{x}{y}$
be the distance from $x$ 
to $y$ in $\T$.  
For the complete graph $K_m$, the distance
$\dT{x}{y}$ is on the order 
of $\sqrt{m}$ and the distribution satisfies
$$\P[\dT{x}{y}> \lambda 
\sqrt{m}]=\exp[-\lambda^2/2]+o(1).$$ 
Moreover, rescaling by dividing edge lengths by $\sqrt m$ 
results in a scaling limit for the UST on 
$K_m$ that is the Brownian {\it continuum random tree} (CRT) 
\cite{AldCRT1}.  
(We will recall the construction of the 
Brownian CRT in \S \ref{CRTsubsection}.)  
Pitman \cite{Pitmanconj} conjectured
that the Brownian CRT should also be the scaling limit
of the UST in certain other graphs, and in particular
should be the scaling limit of the UST on the $d$-dimensional
discrete torus $\Z_n^d$ for large $d$.
This conjecture is supported by a recent result
of Benjamini and Kozma \cite{BenKoz03}, who showed
that the expected distance of {\it loop-erased random walk}
(LERW) between two uniformly chosen points
on $\Z_n^d$ is on the order of $n^{d/2}$ for $d\geq 5$.
Pemantle \cite{Pemantle91} proved that for any $x$ and $y$,
$\dT{x}{y}$
is equal in distribution to the length of a loop-erased random
walk from $x$ to $y$, so Benjamini and Kozma's result also gives some 
information about $\T$.
In Theorem \ref{torusthm}, we confirm that Pitman's conjecture
holds for $d\geq 5$.

\begin{theorem}
\label{torusthm}
Let $x$ and $y$ be independently and uniformly chosen in $\Z_n^d$, 
and let $\dT{x}{y}$ denote the distance
between $x$ and $y$ in the UST. 
For $d\geq 5$,
there exists a constant $\beta(d)>0$ 
such that
\begin{equation}
\label{torustwopointeqn}
\limn \P[\dT{x}{y} > \beta(d) \lambda n^{d/2} ] = 
\exp\left[-\frac{\lambda^2}{2}\right].
\end{equation}
Moreover, if $\{x_1, \dots, x_k\}$ are uniformly chosen from $\Z_n^d$, 
then as $n$ tends to infinity, the joint distribution of 
$$\frac{\dT{x_i}{x_j}}{\beta(d) n^{d/2}}$$
converges in distribution to $F_k$, the joint distribution of distances 
between the first $k$ points of the Poisson line breaking construction of 
the Brownian CRT.
\end{theorem}
In Section \ref{torussection}, we express $\beta(d)$ in the 
form $\gamma(d)/\sqrt{\alpha(d)}$, where $\gamma(d)$ and $\alpha(d)$ are 
probabilities of events involving random walk on $\Z^d$
with $\alpha(d)\rightarrow 1$ and $\gamma(d)\rightarrow 1$ as 
$d\rightarrow \infty$.

Using 
Pemantle's \cite{Pemantle91} result that $\dT{x}{y}$
is equal in distribution to the length of a loop-erased random 
walk from $x$ to $y$,
(\ref{torustwopointeqn}) gives the
limiting distribution of the length of 
LERW from $x$ to $y$ in $\Z_n^d$.
In Subsection \ref{LERWsection},
we will recall Wilson's algorithm
for constructing the UST, which gives
even stronger connections between the UST and LERW.

\subsection{Construction of the CRT and $F_k$}
\label{CRTsubsection}
We define the distribution $F_k$ referred to in Theorem \ref{torusthm},  
via the {\em Poisson line breaking construction\/} of the Brownian CRT.   
Let $s_1, s_2, \dots$ be the arrival 
times for an inhomogeneous Poisson process whose arrival rate at time $t$ 
is $t$.  Draw an initial segment of length $s_1$, and 
label its ends $y_1$ and $y_2$.  Pick a point uniformly on this segment, 
attach a new segment of length $s_2-s_1$, and label the end of this 
segment $y_3$.  
To continue inductively, given a tree with $k$ ends $y_1, y_2, \dots, y_k$
and total edge length $s_{k-1}$, pick a point uniformly on the tree
and attach to this point a
segment of length $s_{k}-s_{k-1}$, and label the new end $y_{k+1}$.  We
let
$F_k$ denote the joint distribution of the distances between the $k$
points $y_1, \dots, y_k$.  Note that $F_2$ is simply
the distribution of $s_1$, which is given by
$$\P[d(y_1, y_2)> \lambda]=\exp\left[-\frac{\lambda^2}{2}\right].$$

As we are only interested in $k$-point distributions of the UST on
graphs, we can stop our construction at time $s_{k-1}$.  More
generally, as $t$ tends
to infinity, the resulting sequence of trees,
viewed as a sequence of metric spaces, converges
to a random, compact metric
space.  The limit is known as the Brownian CRT.  As shown in
\cite{AldCRT1}, for $y_1, \dots, y_k \in K_m$, the joint distribution of
$\dT{y_i}{y_j}/\sqrt m$ converges to $F_k$ as $m\rightarrow \infty$.  Our 
main arguments involve
coupling LERW on $G_n$ with LERW on $K_m$, and we will obtain the results
about the scaling limit of UST on $G_n$ via this coupling.

For a further discussion of the Brownian CRT, see the original
papers of Aldous \cite{AldCRT1, AldCRT2, AldCRT3} or the lecture notes of
Pitman \cite{Pit02}.

\subsection{Other Graphs}

The limiting behavior in (\ref{torustwopointeqn})
is not universal.  On a cycle, for instance, 
the typical distance is on the order of $\vol$ rather than 
$\sqrtvol$.  The case of $\Z_n^2$ is quite hard, but recently Kenyon
\cite{Kenyon00a, Kenyon00b} showed that the typical distance of LERW is on 
the order of $\vol^{5/8}$, and Schramm \cite{Schramm00} and Lawler, 
Schramm, and Werner \cite{LSW} have studied the scaling limit of the UST.

There are also examples of graphs that are not vertex 
transitive,
such as a star, in which the typical distance can be much less than 
$\sqrtvol$. 
Nevertheless, our methods also apply to a 
broader class of graphs, including the hypercubes $\Z_2^n$ and
expander graphs.  As a generalization of Theorem \ref{torusthm},
we will give a set of three conditions such that whenever
all three hold, a suitable scaling limit is again the Brownian CRT.
The three assumptions that we will make are vertex transitivity,
a bounded number of local intersections of two independent random walks, 
and relatively fast mixing of the random walk.  

More formally, let $p^t(x)$ denote the distribution
of simple random walk at time $t$ started at
a basepoint $o \in G_n$.
Our first assumption is that there exists a 
constant $\cc$ such that
\begin{equation}
\label{localstrongtran}
\sup_n \sup_{x\in G_n} \sum_{t=0}^{\sqrtvol} 
(t+1) p^{t}(x) \leq \cc<\infty.
\end{equation}
On $\Z^d$, $p^t(0)$ decays like
$t^{-d/2}$, so for $d\geq 5$, a condition similar to
(\ref{localstrongtran}) holds
in that $\sum_{t=0}^\infty t p^t(0)$ is bounded.  The implication on 
$\Z^d$ is
that two random walks with the same starting point
intersect each other finitely often (see e.g., \cite{Lawbook}, Theorem 
3.5.1).
  Lemma \ref{cutpointlem1}
will show that (\ref{localstrongtran})
is an analog for finite graphs
that says that two independent random walks starting at 
the same point only intersect a bounded number of times in the first 
$\sqrtvol$ steps.

Let $\pi=\pi_n(\cdot)$ denote the stationary distribution
of the walk on $G_n$.
Denote the (uniform) {\it mixing time} by
\begin{equation*}
\mix=
\mixn=\inf\left\{t: \sup_{x\in G_n} 
\left|\frac{p^t(x)}{\pi(x)}-1\right| \leq \frac{1}{2}
\right\}.
\end{equation*}
Our second assumption is that for some $\delta>0$,
\begin{equation}
\label{fastmixing}
\mixn=o(\delsqrtvol).
\end{equation}

Note that this mixing time assumption implicitly requires that the walk be 
aperiodic, as otherwise no such $\mix$ exists.  As adding a holding 
probability of $1/2$ to a random walk does not affect LERW, this 
aperiodicity assumption does not affect the final LERW or the UST.
We will show in Sections \ref{torussection} and \ref{gensection} that
examples of graphs satisfying (\ref{localstrongtran}) 
and (\ref{fastmixing}) include the tori $\Z_n^d$ for $d\geq 5$, the 
hypercube $\Z_2^n$, expanders, and the complete graph $K_n$. 
For vertex transitive graphs, conditions 
(\ref{localstrongtran}) and 
(\ref{fastmixing}) 
are sufficient to generalize \ref{torusthm}.

\begin{theorem}
\label{mainthm}
Suppose that $\{G_n\}$ is a sequence of vertex transitive graphs 
satisfying (\ref{localstrongtran}) and (\ref{fastmixing}) and such that 
$|G_n|\rightarrow \infty$.  Then there 
exists a sequence of constants $\{\beta_n\}$ 
with 
$$
0 < \inf \beta_n \leq \sup \beta_n < \infty
$$
such that
for $x$ and $y$ uniformly chosen from $G_n$,
\begin{equation}
\label{mainboundeqn}
\limn \P[\dT{x}{y} > \beta_n \lambda \sqrtvol ]
= \exp\left[-\frac{\lambda^2}{2}\right].
\end{equation}

Moreover, if $k$ points $\{x_1, \dots, x_k\}$ are chosen independently 
and uniformly from 
$G_n$, then the joint distribution of 
\begin{equation}
\label{CRTeqn}
\frac{\dT{x_i}{x_j}}{\beta_n \sqrtvol}   
\end{equation}
converges to $F_k$ as $n \to \infty$.
 \end{theorem}

Note that the difference between Theorem \ref{torusthm} and Theorem 
\ref{mainthm} is that we potentially need a bounded sequence 
$\beta_n$, rather than having a single scaling constant $\beta$.
If we strengthen assumption (\ref{localstrongtran}), then the factor of 
$\beta_n$ in the rescaling turns out to be identically 1.  More precisely, 
(\ref{localstrongtran}) says that the number of intersections of two walks 
run for $\sqrtvol$ steps is bounded.  To strengthen that, let $\{X_t\}$ 
and $\{Y_u\}$ be independent random walks on $G_n$, 
take $q=(\mix \sqrtvol)^{1/2}$, and suppose that
\begin{equation}
\label{superstrongtran}
\sup_{x\neq y} \E[|\{X_t\}_{t=0}^q \cap \{Y_u\}_{u=0}^q | \mid X_0=x, 
Y_0=y]\rightarrow 0
\end{equation}
as $n\rightarrow \infty$.
Note that condition
(\ref{superstrongtran}) does not hold if the graphs
$\{G_n\}$ have uniformly bounded degree as then the
probability that two walks started at the same point will intersect
after 1 step is bounded away from 0.
It does hold, however, for the hypercubes or the complete graphs.
The local intersections 
are what caused us to require that the points $\{x_i\}$ be chosen 
uniformly, 
and that $\beta_n\neq 1$.  Thus assumption 
(\ref{superstrongtran}) yields the following, stronger 
result.

\begin{theorem}
\label{generalCRTthm}
Let $\{G_n\}$ be a sequence of vertex transitive graphs such that 
(\ref{fastmixing}) and (\ref{superstrongtran}) hold. 
Then for
any choices of $k$ 
distinct vertices $\{x_1^{(n)}, \dots,
x_k^{(n)}\}$ in $G_n$, the joint distributions of
\begin{equation}
\label{generalCRTeqn}
\frac{\dT{x_i^{(n)}} {x_j^{(n)}}}{\sqrtvol}
\end{equation}
converge to $F_k$ when $n \to \infty$.
\end{theorem}

In Section \ref{LERWsection} we define LERW and recall Wilson's 
algorithm.  In Section \ref{newsketch} we give an outline
of the proof of Theorem \ref{torusthm}, explain the
significance of the size $\sqrtvol$, and discuss the structure   
of the rest of the paper. 

\subsection{Loop-erased random walks and Wilson's algorithm}
\label{LERWsection}

Here and throughout this paper, we will use $\langle \cdot \rangle$ to 
denote sequences when order is important, and $\{ \cdot \}$ to denote 
sets.

Given a finite path $\gamma=\langle u_0, u_1, \dots, u_\ell \rangle$ in a 
graph  
$G$, let $\LE(\gamma)$ denote the {\it loop-erasure} of $\gamma$ with
loops erased in chronological order.  Formally,
$\LE(\gamma)$ is the sequence $\langle v_0, v_1, \dots \rangle$
constructed recursively as follows:
first, let $v_0=u_0$; then, given
$v_r$, let $k$ be the last time that $u_k=v_r$,
and let $v_{r+1}=u_{k+1}$, with the convention that if $k=\ell$,
then $v_r$ is the last term of $\LE(\gamma)$.
In the case when $\gamma$ is the path of a random walk that starts at 
$x$ and is stopped when it reaches $y$, we call $\LE(\gamma)$ a 
{\it loop-erased random walk} (LERW) from $x$ to $y$.  More generally, if 
$\gamma$ is a random walk from $x$ stopped when it hits a set $S$, then we 
say $\LE(\gamma)$ is LERW from $x$ to $S$.

For an unweighted graph $G$, the UST on $G$ is connected to 
LERW on $G$ by Wilson's algorithm for constructing a UST: pick a root 
vertex 
$\rho$, and form an initial subtree $\T_1$ by picking another vertex $x_1$ 
and running LERW from $x_1$ to $\rho$.  We then proceed recursively as 
follows: given a subtree $\T_i$, pick a vertex $x_{i+1}$ and run LERW from 
$x_{i+1}$ to $\T_i$.  Let $\T_{i+1}$ be the union of $\T_i$ with this new 
path.  Proceeding until $\T_i$ is a spanning tree yields a UST 
\cite{Wilson}. 
Taking $\rho=y$ and $x_1=x$ yields Pemantle's result
\cite{Pemantle91} that the 
distribution of the path from $x$ to $y$ within the UST is the same as 
the distribution of LERW from $x$ to $y$.
We will use the fact that
Wilson's algorithm is very robust--in particular, the sequence $\{x_i\}$
may be chosen arbitrarily, and it also applies to weighted graphs.
To apply Wilson's algorithm to weighted graphs, we let LERW denote the 
loop-erasure of a weighted random walk that moves according to the weights 
on the graphs.  The probability of any given spanning tree resulting from 
Wilson's algorithm then is proportional to the product of the weights 
of the edges of the tree, so Wilson's algorithm yields a weighted UST
\cite{Wilson}.

As mentioned earlier, the limiting distribution of $\dT{x}{y}$ on 
the complete graph $K_m$ is given by
$$\P[\dT{x}{y}> \lambda 
\sqrt{m}]=\exp[-\lambda^2/2]+o(1).$$ 
Because there is no geometry on the complete graph, 
there are a number of derivations of this limit
(c.f., \cite{CamPit}).
One such argument that uses the connection between the UST and
the LERW is as follows:
Let $\gamma=\langle x_1, \dots, x_\ell \rangle$ denote the loop-erased 
path from $x$ to $y$, so $x_1=x$ and $x_\ell=y$.
Conditioned on $\gamma_i=\langle x_1, \dots, x_i \rangle$, we want to 
compute $\P[x_{i+1}=z]$.  Considering the walk after the last visit to 
$\gamma_i$, we want to condition on never returning to $\gamma_i$.
Let $f(z)$ be the probability that the next step is to $z$ 
and the walk then reaches $y$ before hitting $\gamma_i$.  
We have $f(y)=1/(m-1)$, and
by symmetry,
$f(z)=1/[(m-1)(i+1)]$ for $z\in K_m \setminus\{\gamma_i \cup y\}$.  As 
$f\equiv 0$ on $\gamma_i$, conditioning on the union of these events 
occurring show that
the probability of stepping to 
$y$ at step $(i+1)$ is thus $(i+1)/m$.  In particular,
$$
\P[\dT{x}{y}=k] = \frac{k+1}{m} 
\left[\prod_{i=1}^{k-1}\left(1-\frac{i+1}{m}\right) \right].
$$
Rescaling by $m^{-1/2}$  and letting $m\rightarrow \infty$,
this distribution converges to the Rayleigh distribution.
Generalizing this argument for graphs with more structure
is difficult, so we will use other techniques.

\section{Outline of Proof}
\label{newsketch}

To prove Theorem \ref{mainthm}, we would like to use Wilson's 
algorithm and make a comparison with the Poisson line breaking 
construction.  Unfortunately, if $x$ and $y$ are 
uniformly chosen from a graph $G_n$ that satisfies the hypotheses 
of Theorem \ref{mainthm}, it turns out that the typical random walk from 
$x$ to $y$ takes on the order of $\vol$ steps.
In the special case of $G_n=\Z_n^d$, Benjamini and Kozma \cite{BenKoz03}
showed that the typical length of the loop-erased path is only of the 
order of $\sqrtvol$, meaning that almost all of the original path is 
erased.
Theorem \ref{mainthm} implies that this problem occurs
more generally, so to avoid it we would prefer to only run a 
random walk for the 
order 
of $\sqrtvol$ steps.

To do this, for a given $L$, consider the extension $\extG{n}{L}$ 
of $G_n$ 
formed by adding a vertex $\rho$, with an edge from every vertex of $G_n$
to $\rho$ of weight such that a weighted random walk on $\extG{n}{L}$ is
simple
random walk on $G_n$ modified to move to $\rho$ after a geometric number
of steps with mean $L \sqrtvol$.  We 
will run Wilson's algorithm on $\extG{n}{L}$ with root vertex $\rho$.  
This 
results in a spanning tree $\extT$ on $\extG{n}{L}$, rather than the 
uniform 
spanning tree $\T$ on $G_n$.  The tree $\extT$ induces a spanning forest 
$\F$ on $G_n$ formed by restricting $\extT$ to $G_n$.  
Let $\dT{x}{y}$ and $\dF{x}{y}$ denote the distance between $x$ 
and $y$ in $\T$ and $\F$ respectively, with the convention that 
$\dF{x}{y}=\infty$ if $x$ and $y$ are in different components of $\F$.
This spanning 
forest $\F$ is comparable in the following sense to the UST $\T$:

\begin{lemma}
\label{mainstochasticdomlem}
Let $\{G_n\}$ be a sequence of vertex transitive graphs 
such that (\ref{localstrongtran}) and (\ref{fastmixing}) hold with 
constants $\cc$ and $\delta$ respectively, and such that $|G_n|\rightarrow 
\infty$.
Pick $k$ points $\{x_1, \dots, x_k\}$ uniformly and independently from 
$G_n$.
For any $\eps>0$, 
there exist constants $L_0$ and $N$ such that 
if 
$L>L_0$ and $n>N$, then the total variation distance between the joint 
distribution of the distances $\dT{x_i}{x_j}$ and the joint 
distribution of the distances $\dF{x_i}{x_j}$ is less than $\eps$.
\end{lemma} 

We will prove Lemma \ref{mainstochasticdomlem} in Section 
\ref{treesection}.  The point of Lemma \ref{mainstochasticdomlem}
is that it now suffices to understand $\F$ rather than $\T$.
To do so, we will use a comparison 
 with an extension of the complete graph $K_m$.
Let $K_m$ denote the complete graph with a self-loop at every vertex,
and construct an extension
$\extK{m}{\widetilde L}$ of $K_m$ in the same way that 
$\extG{n}{L}$ was 
constructed, meaning that $\extK{m}{\widetilde L}$ is
formed by adding a vertex $\rho$ that is connected
to every vertex of $K_m$ with an edge of
weight $m/(\widetilde L \sqrt{m}-1)$ and so
weighted random walk on $\extK{m}{\widetilde L}$ is simple random
walk on $K_m$ that jumps to $\rho$ after a geometric number
of steps with mean $\widetilde L \sqrt{m}$.
  
On $G_n$, pick $\{x_1, 
\dots, x_k\}$ uniformly and independently.  Running Wilson's 
algorithm on $\extG{n}{L}$ with root vertex $\rho$,
let $\T_0=\{\rho\}$.
Given $\T_\ell$, let $\T_{\ell+1}$ be formed by running LERW from 
$x_{\ell+1}$ to $\T_\ell$.
Note that $\T_1$ is LERW on $\extG{n}{L}$ from $x_1$ to $\rho$.  
We call $\T_\ell$ the $\ell$-th partial spanning tree.
Likewise, let
$\widetilde \T_\ell$ denote the $\ell$-th partial spanning tree on 
$\extK{m}{\widetilde L}$.
We want to select $m$ and $\widetilde L$ in such a way to couple
$\T_\ell$ and $\widetilde \T_\ell$.
To couple both hitting probabilities and length, we need to
consider the following definition of {\it capacity}.

\begin{definition}
\label{capdef}
Let $\{X_t\}$ be simple random walk on $G_n$.  For $S\subset G_n$,
let $T_S=\min\{t\geq 0: X_t\in S\}$ denote the hitting time of $S$.
The $r$-capacity of $S$ is given by
$$
\Cap(S)=\P_{\pi} [T_S < r].
$$
\end{definition}

By considering the expected number of visits
to $S$ in the interval $[0,r)$,
\begin{equation}
\label{easycapbound}
\Cap(S) \leq r \pi(S) = \frac{r |S|}{\vol}
\end{equation}
 for any set $S$.  When $S$ is a segment of a random walk, we 
will show in Lemma \ref{concentrationlem} 
that the bound in 
(\ref{easycapbound}) is, with high probability, sharp up to constants.

\begin{definition}
\label{couplingdef}
Let $\{x_1, \dots, x_k\}$ and $\{y_1, \dots, y_k\}$ be
uniformly and independently chosen from $G_n$ and $K_m$
respectively.
For any $C$, $L$, $\widetilde L$, $\alpha$, $\beta$, $\delta$, and $r$,
the partial spanning trees $\T_k$ and $\widetilde \T_k$ on $\extG{n}{L}$ 
and $\extK{m}{\widetilde L}$ are said to be successfully 
coupled if
the following three conditions hold:\\
\begin{equation} 
\label{couplecap}
\left| \frac{\Cap(\T_k) \sqrtvol}{\alpha^{1/2} r} - \frac{|\widetilde
\T_k|}{\sqrt{m}} \right| \leq
C \vol^{-\delta/32}
\end{equation}
for all $i,j \leq k$, 
\begin{equation}
\label{couplelength}
\left| \frac{\dTk{x_i}{x_j}}{\beta \sqrtvol} 
- \frac{\dtildeTk{y_i}{y_j}}{\sqrt m} \right| \leq C \vol^{-\delta/32} 
\end{equation}
Let $\mu$ be the uniform measure on $\T_k$ and $\nu$ 
the hitting measure, 
meaning that if $\{X_t\}$ is simple
random walk on $G_n$ and $T=\inf\{t\geq 0: X_t\in \T_k\}$,
then $\nu(x)=\P[X_T=x]$.
There is a coupling of $\mu$ and $\nu$ such that if
$\xi$ is chosen according to $\mu$ and $\eta$ according to $\nu$, then
 \begin{equation}
\label{uniformhitting}
\P[\dTk{\xi}{\eta}< C \vol^{1/2-\delta/32}]\geq 1-C(\vol^{-\delta/32}).
\end{equation}
\end{definition}

\begin{lemma}
\label{couplinglemma}
Let $\{G_n\}$ be a sequence of vertex transitive graphs with 
$|G_n|\rightarrow \infty$
such that (\ref{localstrongtran}) and (\ref{fastmixing}) hold with 
constants $\cc$ and $\delta$ respectively.
Then there are sequences of constants $\alpha=\alpha(n)$,
$\beta=\beta(n)$ and $r=r(n)$, and
for any $L$ and $\eps>0$, there are sequences 
$\widetilde L=\widetilde L(L,n)$ and 
$m=m(n)$, and a constant $C=C(k)$
such that if
$\{x_1, \dots, x_k\}$ and $\{y_1, \dots, y_k\}$ are chosen 
independently without replacement from $G_n$ and $K_m$ respectively,
then for any $n>N$, there is a coupling of $\T_k$ and $\widetilde 
\T_k$
that is successful with probability 
$1-C(k) \left(\vol^{-\delta/32}\right)$.
Moreover, $m$ and $\widetilde L$ can be made arbitrarily large by
first taking $L$ and then $n$ sufficiently large.
 \end{lemma}

Theorem \ref{mainthm} is then a consequence of (\ref{couplelength}),
along with Lemma \ref{mainstochasticdomlem} and Aldous' theorem
that the Brownian CRT is the scaling limit of UST on $K_m$.
Although only (\ref{couplelength}) is needed to imply Theorem
\ref{mainthm}, we add conditions (\ref{couplecap}) and 
(\ref{uniformhitting}) to Lemma \ref{couplinglemma}
because the proof will be by induction, and the
inductive step requires those two pieces.

The key to the proof
is a rescaling argument that makes use of the facts that loops formed by a 
typical random walk are either very short or very long.  More formally, 
let $\{X_t\}$ be a simple random walk on $G_n$ run for $L \sqrtvol$ steps.  
Call a loop short if it is of length at most $\mix$, and long if it is of 
length at least $\vol^{1/2-\delta}$.  If $t-u \geq \mix$, then
$\P[X_t=X_u] \leq 2/\vol$, so the expected number of loops of length 
between $\mix$ and $\vol^{1/2-\delta}$ formed by time $L\sqrtvol$ is 
bounded by $2 L \vol^{-\delta}$.  In particular, the probability of having 
loops of this type of intermediate length tends to 0.

 In Section \ref{sketch} we will describe a rescaling argument that takes 
advantage of the absence of intermediate length loops.  In 
order to make this rescaling argument work, we will show in Section 
\ref{lemmasection} that short pieces of a loop-erased random walk retain a 
positive proportion of their length.  Section
\ref{largedevsection} then combines these pieces of the walk to show
that longer lengths of loop-erased walk have a behavior that is tightly
concentrated around the mean behavior.  
Section \ref{couplingsection} is then devoted to using these pieces to 
prove Lemma \ref{couplinglemma}.
In Section \ref{treesection}, we 
prove Lemma \ref{mainstochasticdomlem}, thus concluding
the proof of Theorem \ref{mainthm}.
In Section \ref{torussection} we then prove Theorem \ref{torusthm},
and in Section \ref{gensection} we prove Theorem \ref{generalCRTthm}.

\section{Key definitions and introduction of rescaling}
\label{sketch}

As the inductive step in the proof of Lemma \ref{couplinglemma} requires 
running LERW from $x_{k+1}$ to $\T_k$ and showing that this new 
length is similar to what occurs on $\extK{m}{\widetilde L}$, we 
begin by 
considering what happens on $\extK{m}{\widetilde L}$.

Suppose that $S\subset \extK{m}{\widetilde L}$ is a subset such 
that $\rho \in S$.  Pick a point $y\in \extK{m}{\widetilde L}$
with $y\neq \rho$, and let $\{Y_u\}$ 
denote 
weighted random walk on $\extK{m}{\widetilde L}$ started at $y$ and 
stopped at 
time 
$$T=\min\{u \geq 0: Y_u\in S\}.$$
We wish to understand
$
\LE \langle Y_u\rangle_{u=0}^{T},
$
and in particular its length.
To do this, let
\begin{equation}
\label{tildeIijdef}
\tildeIij=\begin{cases}
1 & Y_i=Y_j \\
0 & Y_i\neq Y_j \\
\end{cases}
\end{equation}
be a collection of indicator random variables that keep track of 
intersections.
One special property of the complete graph is that,
for fixed $j$, the joint distribution 
of $\tildeIij$ as $i$ varies
conditioned on the values of $\tildeIkl$ for $k, \ell <j$
is the same as the joint distribution conditioned on 
$\{Y_u\}_{u<j}$.
Since we are only interested in the length of LERW from $y$ to $S$,
the lack of geometry on $K_m$ means that we 
do not 
lose information by only keeping
track of the time indices instead of locations of $\{Y_u\}$.  Moreover,
keeping track of times instead 
of locations will generalize more 
appropriately for our rescaling argument when
 we consider LERW on $\extG{n}{L}$.

We now inductively construct a family of sequences 
$\{\Kindex{j}\}$ that record what time indices have survived 
loop-erasure up to time $j$.  
 Let our initial sequence of length one be given by
$
\Kindex{0}=\langle 0\rangle.
$
Letting $\concat$ denote concatenation of sequences,
for $0<j\leq T$ let
\begin{equation}
\label{Kindexdef}
\Kindex{j}=\langle k \in \Kindex{j-1}: \tildeIij=0 \ \forall i 
\leq k, 
i \in \Kindex{j-1} \rangle \concat \{j\}.
\end{equation}
The result of this definition is that
 $\Kindex{j}$ consists of the original time
indices of the walk that have survived loop-erasure at time $j$, 
with the convention that when a loop is formed, the original time index is 
removed and the new time is retained.  
One consequence of this is that 
for $j\leq T$,
$$
\LE \langle Y_u \rangle_{u=0}^j = \langle Y_u \rangle_{u\in 
\Kindex{j}}.
$$

For the random walk on $\extG{n}{L}$, we will again define a collection of 
indicator random variables $\{\Iij\}$ and index sets $\{\Gindex{j}\}$, 
but the difference now is that the indices will represent a moderately 
long 
segment of the loop-erased random walk instead of individual points.

In the rest of this paper, we will be using a variety of different time 
scales;
for ease of reference we give here a summary of the 
meanings of the different
scales on which we will be working.  First,
weighted random walk on the extension $\extG{n}{L}$
started in $G_n$ moves to $\rho$
after a geometric number of steps that is on the order 
of $\sqrtvol$.  Recall that $\mix=\mixn$ denotes the mixing time for 
random walk on $G_n$. 
We will break up $\sqrtvol$ into shorter segments $A_i$ of length 
that is asymptotically
$r=\lfloor \mix^{1/4}\vol^{3/8}\rfloor$.
By assumption (\ref{fastmixing}), $\mix$, and thus $r$, is of a lower 
order than $\sqrtvol$.
To show that the behavior of the walk on
each run of length $r$ is close to its mean
behavior, we will further break these runs into smaller
pieces of length $q=\lfloor \mix^{1/2}\vol^{1/4}\rfloor$ and then sum the 
pieces to get
large deviation estimates.  These estimates
require independence between the segments, so instead
of considering the loop-erased path, we will consider
a local loop-erasure, with a window size $s=\lfloor \mix^{3/4}\vol^{1/8} 
\rfloor$ that is much smaller than $q$.  To justify that the
restriction to local loop-erasure does not throw away too much 
information, we will then finally
show that there are a number of local cutpoints,
for which we will only look at path segments whose length is on the order 
of $\mix$.
As a summary, see Table \ref{deftable}.

\begin{table}[h,t]
\begin{center}
\begin{tabular}{|c|c|} \hline
$\sqrtvol$ & typical length of LERW \\ \hline
$r=\lfloor \mix^{1/4} \vol^{3/8} \rfloor
$ & length of $A_i$ \\ \hline
$q = \lfloor \mix^{1/2} \vol^{1/4} \rfloor
$ & subdivisions of $A_i$ \\ \hline
$s = \lfloor \mix^{3/4} \vol^{1/8} \rfloor$ 
& window for local loop-erasure \\ \hline
$\mix=\mixn$ & mixing time, and window for local cutpoints \\ \hline
\end{tabular}
\end{center}
\caption{Definitions of scales}
\label{deftable}
\end{table}

To remember these relative sizes, note that 
$\{\mix, s, q, r, \sqrtvol\}$ is approximately a
 geometric sequence with common ratio 
$\mix^{-1/4} \vol^{1/8}$.

As mentioned above, in order to get the independence
between various pieces that we need for our
large deviation estimates, we will need to use
a local loop-erasure instead of the original loop-erasure.
\begin{definition}
\label{localretaineddef}
Let $\{X_t\}$ denote a weighted random walk on $\extG{n}{L}$ with $X_0\neq 
\rho$ that is stopped at time $T$.
A time $u$ is {\bf locally retained} if
$$\LE \langle X_t \rangle_{t=\max\{0,u-s\}}^u \cap
\{X_t\}_{t=u+1}^{\min\{T,u+s\}} = \emptyset.$$
\end{definition}
\begin{definition}
\label{LLEdef}
For a stopping time $T$,
let $U$ be the set of all times $t\leq T$ that are locally retained.  The 
{\bf local loop-erasure} $\LLE \langle X_t \rangle_{t=0}^T$ is the 
subsequence 
of $\langle X_t \rangle_{t=0}^T$ such that
$$
\LLE \langle X_t \rangle_{t=0}^T := \langle X_u \rangle_{u\in U}.
$$
\end{definition}
When the stopping time $T$ is not specified, we will take $T=\infty$
in the definition of the local loop-erasure.

It is not {\it a priori} true that either the local loop-erasure
or the loop-erasure contains the other.
For example, if the original path has a loop of length
$s-1$, the local loop-erasure could have a jump,
while if the original path has a loop of length just longer
than $s$, then the local loop-erasure can have short loops.
These differences raise a problem that we will 
have to deal with later.
We will later formalize the notion that, with high probability, 
the local loop-erasure has no jumps, and that the main 
difference between the local loop-erasure and the loop-erasure comes from 
having long loops (meaning of length greater than $r$).  Our 
coupling with LERW on the complete graph will keep track of the long 
loops.

Given a set $S\subset \extG{n}{L}$ with $\rho \in S$,
pick a point $x\in G_{n}$, and let $\{X_t\}$ be weighted random walk 
on 
$\extG{n}{L}$ with $X_0=x$ and run until time $T=\min\{t\geq 0: X_t\in 
S\}$.

For $i< T/r$, let
\begin{equation}
\label{Aidef}
A_i=A_i(r,s)=\Z \cap [(i-1)r+2s+1, ir-s].
\end{equation}
As $s=o(r)$, the number of elements of $A_i$ is asymptotically $r$.
Adding a buffer of length $s$ 
at the beginning and end of $A_i$ means that
the times that are locally retained within the different $A_i$
are independent.  The second delay of $s$ at the start of $A_i$
will also mean that the locations of the path on different $A_i$
are close to independent.
We will let $\LLE(A_i)$ denote
the part of the local loop-erasure of $\{X_t\}$ whose original times
were in $A_i$, that is to say
$$\LLE(A_i) := \langle X_t \rangle_{t\in A_i \cap U}.$$

Since $\rho\in S$, the hitting time
$T$ is on the order of at most $\sqrtvol$, so the number of intervals 
$A_i$ is at most on the order of 
$\sqrtvol/r$.
We now wish to keep track of non-local loops, meaning loops
that somehow involve two of the $A_i$. 
For $i<j$, let $\Iij$ denote indicator random variables
\begin{equation}
\label{Iijdef}
\Iij=\begin{cases}
1 & \{\LLE (A_i) \cap
\{X_t\}_{t\in A_j}\neq \emptyset\}, \\
0 & \{\LLE (A_i) \cap
\{X_t\}_{t\in A_j}= \emptyset\}. \\
\end{cases}
\end{equation}
Unlike the case of the complete graph, here
the joint distribution of $\Iij$ conditioned
on $\{X_t\}_{t < rj}$ is different
from the joint distribution of $\Iij$ conditioned
on $\{\Ikl, k,\ell<j\}$. 
 Despite this, we can still recursively
construct a family of sequences $\Gindex{j}$ that in some
sense records which runs $\LLE(A_i)$ survive
loop-erasure.
The construction will implicitly take
the entire path $\{X_t\}$ into consideration.  To begin, let 
$\Gindex{0}=\langle 0 \rangle$.
For $0<j \leq \lceil T/r \rceil$,
let $\concat$ denote concatenation of sequences, and let
\begin{equation}
\label{Gindexdef}
\Gindex{j}=\langle k \in \Gindex{j-1}: \Iij=0 \ \forall i
\leq k, i \in \Gindex{j-1} \rangle \concat \{j\}.
\end{equation}

These sequences are intended to play much the same role as 
$\{\Kindex{j}\}$, 
and can be thought of as keeping track of indices $i$ such that 
$\LLE(A_i)$ 
is completely
contained inside $\LE \langle X_t \rangle_{t=0}^T$.  There are some 
problems:
if $\LLE(A_i)$ is involved in a long loop, then part of it is erased.
In particular,
it is not true that $\LLE(A_i) \subset \LE \langle X_t 
\rangle_{t=0}^{T}$, and also if the long loop involves $\LLE(A_j)$,
then only one of $i$ or $j$ is retained in $\Gindex{j}$.
There is also a problem when one end of a long loop is in the
gap of length $3s$ between the $A_i$.
We will later prove that these differences are sufficiently rare that 
their
contribution is of a lower order of magnitude than the length of LERW.

In Section \ref{largedevsection}, we will show that the 
length of $\LLE(A_i)$ is tightly concentrated about its mean $\gamma 
r$ for some $\gamma=\gamma(n)$ that is bounded away from $0$.  
This implies that 
$|\LE \langle X_t \rangle_{t=0}^T|$
is roughly $\gamma r |\Gindex{\lfloor T/r \rfloor}|$.
We will prove the induction step of Lemma \ref{couplinglemma}
in Section \ref{couplingsection}, but to understand the
idea, take $S=\T_k$.
We want to pick an $m$ that lets us couple LERW on $\extG{n}{L}$ and 
$\extK{m}{\widetilde L}$ in 
such a way that $|\Gindex{\lfloor T/r\rfloor}|=|\Kindex{\widetilde T}|$,
where $\widetilde T$ is the hitting time for  
$\widetilde \T_k \subset \extK{m}{\widetilde L}$.
To do this, 
we need the distribution of $\Iij$ to be close to the
distribution of $\tildeIij$.
Let $\Cap(S)$ denote the capacity of $S$ as defined
in Definition \ref{capdef}.

Define $\gamma=\gamma(n)$ and $\alpha=\alpha(n)$ by 
\begin{equation}
\label{alphagammadef}
\gamma=\E \frac {|\LLE(A_i)|}{r}, \qquad
\alpha=\E \Cap[\LLE(A_i)]\frac{\vol}{r^2}.
\end{equation}
In Lemma \ref{concentrationlem}, 
we will show that both $|\LLE(A_i)|$ and $\Cap[\LLE(A_i)]$ are tightly 
concentrated about their means.  
The capacity bound implies that
$\E \Iij$ is approximately $\alpha r^2/\vol$.  On $K_m$,
$\E \tildeIij=(1+o(1))/m$, 
so in our rescaling we will take the size of the 
complete graph to be
$$
m=\left\lceil \frac{\vol}{\alpha r^2}\right\rceil.
$$
and let
$$
\beta=\beta_n=\frac{\gamma}{\alpha^{1/2}}.
$$
For the special case of $S=\T_0=\{\rho\}$, taking $\widetilde T$
to be the hitting time of $\rho$ on $\extK{m}{\widetilde L}$,
the two point distribution reduces to showing that
\begin{align*}
\P\left[ | \LE\langle X_t \rangle_{t=0}^T | > 
\frac{\gamma}{\alpha^{1/2}} \lambda \sqrtvol\right]
& = \P \left[ |\Gindex{\lfloor T/r\rfloor} |> \frac{
\lambda \sqrtvol}{\alpha^{1/2}r} \right]+o(1) \\
& = \P \left[ | \Kindex{\widetilde T} > \lambda \sqrt{m} \right]
+o(1),
\end{align*}
which tends to $\exp[-\lambda^2/2]$ as $m, \widetilde L \rightarrow
\infty$.

To obtain (\ref{torustwopointeqn}), in Section
\ref{torussection} we will show that, on the torus,
the limits of $\gamma(n)$ and $\alpha(n)$ exist.
Computing these limits
and replacing $\beta_n$
by the limit of $\beta_n$ will then prove Theorem \ref{torusthm}.

In the next few sections, we will develop the tools needed
to prove Lemma \ref{couplinglemma}.

\section{Positive length of small pieces}
\label{lemmasection}

The aim of this section is to study local-loop erasure of runs of 
length roughly
$q\approx \mix^{1/2} \vol^{1/4}$ and show that
 $\LLE \langle X_t \rangle_{t=2s+1}^{q-s}$ retains a positive proportion
of the original walk.  Note that the buffers of length $s$ and $2s$ 
at the start and end of these pieces of length $q$ are the same size as in 
$A_i$ and serve the same role.
The fact that erasing loops shortens the path, along with
(\ref{easycapbound}), gives the upper bounds
\begin{align*}
\frac{\E \left|\LLE \langle X_t \rangle_{t=2s}^{q-s}\right|}{q} & \leq 
1, \\
\E \Cap[\LLE \langle X_t \rangle_{t=2s}^{q-s}] \frac{\vol}{qr}& \leq 1,
\end{align*}
where $\Cap(S)$ is as in Definition \ref{capdef}.   
The focus of this section will be giving lower bounds for these 
quantities, 
and in particular showing that they are bounded 
away from 0.

For random walks on $\Z^d$, the fact that $\sum k \P[X_k=0]<\infty$ 
 for $d\geq 5$ is equivalent to the fact that two simple random 
walks on $\Z^d$ will intersect each other finitely often in dimensions 5 
and higher.  The next lemma makes more precise the fact that condition 
(\ref{localstrongtran}) gives a local analog.

\begin{lemma}
\label{cutpointlem1}
Suppose that $\{G_n\}$ is a sequence of 
vertex transitive graphs 
satisfying (\ref{localstrongtran}) with constant $\cc$.
Let $\{X_t\}$ and $\{Y_t\}$ be independent
random walks started at the same point $o$.  Let $\{\widetilde X_t\}$ and 
$\{\widetilde Y_t\}$ be the walks $\{X_t\}$ and $\{Y_t\}$ killed at 
random times
$T_X$ and $T_Y$, which are geometrically distributed
random variables with mean $(1-\lambda)^{-1}$.
Letting $\cc$ be the constant from (\ref{localstrongtran}),
\begin{equation}
\label{cutpointeq}
\P(\{\widetilde X_t\}_{t\geq 0} \cap \{\widetilde Y_t\}_{t\geq 
1}=\emptyset)
\geq \left(\cc+\frac{2}{\vol(1-\lambda)^2} 
\right)^{-1}.
\end{equation}
\end{lemma}

\begin{proof}
We use the central idea of the proof of Proposition 3.2.2 in 
\cite{Lawbook}.  Call a pair of times $(i,j)$ a *-last intersection if 
$$\{(t,u): \widetilde X_t=\widetilde Y_u, t\geq i, u\geq j\}=\{(i,j)\}.$$
Abbreviate 
$$f(\lambda)=\P(\{\widetilde X_t\}_{t\geq 0} \cap \{\widetilde 
Y_t\}_{t\geq 1}=\emptyset).$$

Note that $\P[\widetilde X_i=\widetilde Y_j]=\lambda^{i+j} 
\P[X_i=Y_j]$, so
by vertex transitivity and the memoryless
property of exponential random variables,
the probability that $(i,j)$ is a *-last intersection is bounded
by
$$
\lambda^{i+j}
\P[X_i=Y_j] f(\lambda).
$$
Because the killed paths are finite in length,
there is at least one *-last intersection.
By symmetry of the walk, $\P[X_i=Y_j]=\P[X_{i+j}=o]$, so
considering the expected number of *-last intersections gives
\begin{align*}
1 & \leq 
\sum_{i=0}^\infty \sum_{j=0}^\infty \lambda^{i+j} \P[X_i=Y_j] f(\lambda) 
\\
& = \sum_{k=0}^\infty \lambda^k (k+1) \P_o[X_k=o] f(\lambda) \\
& \leq \left( \sum_{k=0}^{\sqrtvol} (k+1) \P_o[X_k=o] + 
\sum_{k=\sqrtvol}^\infty \frac{2(k+1)}{\vol} 
\lambda^k \right) f(\lambda) \\
& \leq \left(\cc + \frac{2}{\vol (1-\lambda)^2}\right) f(\lambda).
\end{align*} 
\end{proof}

For a path $\langle X_t \rangle_{t=0}^T$, a point $X_u$ is a {\it local 
cutpoint} if
$\{X_{u-\mix}, \dots, X_{u-1} \} \cap \{X_{u+1}, \dots, X_{u+\mix}
\}=\emptyset$.
Lemma \ref{cutpointlem1} and condition (\ref{fastmixing})
imply that the expected number
of local cutpoints of random walk is a positive proportion of the 
length of the path:   
taking $1-\lambda=(\mix \sqrtvol)^{-1/2}$ means that $[(1-\lambda)^2 
|G_n|]^{-1}$ tends to 
$0$, and also that $\P[\min\{T_X, T_Y\} >\mix]=1-o(1)$.
Lemma \ref{cutpointlem1} then implies that the probability that 
a given point is a local cutpoint is at least $1/\cc + o(1)$.
  The significance of 
local cutpoints is that
conditioned on not having long loops (meaning loops of 
length greater than 
$\mix$), local cutpoints are also retained in the loop-erasure.
\begin{corollary}
\label{cutptcor}
Suppose that $\{G_n\}$ is a sequence of 
graphs satisfying the assumptions of Theorem \ref{mainthm}, 
and let $\{X_t\}$ be simple random 
walk on $G_n$ with geometric killing rate $(L \sqrtvol)^{-1}$.  Let $T$ be 
the killing time, and let $\Gamma_n$ be the event that for all 
$t\in[0,T-s]$, the interval $[t,t+s]$ contains a local cutpoint.
Then there is a constant $C$ such that for any $s$, 
 $\P[\Gamma_n]\geq 1-L \sqrtvol \exp(-C s/\mix)$.
\end{corollary}

Note that taking $s\geq (\log \vol)^2 \mix$, as the choice in Table 
\ref{deftable} does, gives $\P[\Gamma_n]\geq 1-o(\vol^{-1})$.

\begin{proof}
By Lemma \ref{cutpointlem1}, any point is a local cutpoint with 
probability bounded away from 0.  Whether or not $X_t$ and $X_{t+2\mix}$ 
are local cutpoints 
are independent events, so the probability of having no local cutpoint in 
an interval of length $2k\mix$ decays exponentially in $k$.  The result 
then follows from  
summing over all start times $t\in [0, T-s]$.
\end{proof}

Repeating the argument that the the probability of a point being a
local cutpoint is bounded below,
but considering times that are locally retained rather than local 
cutpoints, shows
that the probability that any given time is locally retained
is at least $(1+o(1))/\cc$.  In particular, this gives

\begin{corollary} Suppose that $\{G_n\}$ is a sequence of 
vertex transitive graphs 
satisfying (\ref{localstrongtran}) with constant $\cc$,
and suppose that $s=s(n)$ and $q=q(n)$ are such that $q\leq 
\vol^{1/2-\delta/4}$ and 
$s=o(q)$.  If $\{X_t\}$ is 
simple random walk on $G_n$, then
\label{explengthcor}
\begin{equation}
\label{explengtheqn}
\E\left[ | \LLE\langle X_t \rangle_{t=2s}^{q-s}| \right] 
\geq \frac{q+o(q)}{\cc}.
\end{equation}
\end{corollary}

In addition to knowing that a positive fraction of the walk is retained by 
local loop-erasure, we also want to know that the probability of two 
random walks intersecting is not too greatly reduced by local loop-erasure 
of one of the paths.  
For two i.i.d., transient random walks $\{X_t\}$ and $\{Y_u\}$,
this holds quite generally.  For example,
\begin{equation}
\label{LPSresult}
\P[\LE\langle X_t\rangle \cap \{Y_u\} \neq \emptyset] 
\geq 2^{-8} \P[\langle X_t\rangle \cap \{Y_u\} \neq \emptyset]
\end{equation}
(see \cite{LPS03}).  
Markov chains with geometric killing on a finite state space
are transient Markov chains, but (\ref{LPSresult})
does not apply in our case, partly because we need 
deterministic (rather than geometric) killing, but 
mostly because we are interested in cases when the killing times are on 
different orders of magnitude and so the killed walks are not i.i.d.

\begin{lemma}
\label{intersectionlemma}
Suppose that $\{G_n\}$ is a sequence of 
vertex transitive graphs 
satisfying (\ref{localstrongtran}) with constant $\cc$.
Let $\{X_t\}$ and $\{Y_u\}$ be two independent random walks 
run with deterministic killing times $T_X$ and $T_Y$ respectively, run 
from uniformly chosen starting locations $x, y\in G_n$.
If $T_X, T_Y \leq \sqrtvol/2$
and $s\leq \eps \sqrtvol$, then 
\begin{equation}
\P[\LLE \langle X_t \rangle_{t=0}^{T_X} \cap \{Y_u\}_{u=0}^{T_Y} \neq 
\emptyset] \geq \frac{T_X T_Y}{\vol \cc^2(4\cc-3)} 
\left[\frac{1}{1+2\eps}-2(1-\sqrt{\eps})\right]
\end{equation}
where $\cc$ is the constant in (\ref{localstrongtran}).
\end{lemma}

\begin{proof}
The proof relies on the second moment bound 
$$\P[Z>0]\geq \frac{(\E Z)^2}{ \E(Z^2)}$$
 for non-negative random variables $Z$.
 For $i\leq T_X$ and $j\leq T_Y$, 
let $J_{ij}$ be an indicator random variable for the event 
$$\{X_i=Y_j \in \LLE \langle X_t \rangle_{t=0}^{T_X}\},$$
and take $Z=\sum_{i,j} J_{ij}$ to be the number of ordered pairs $(i,j)$ 
corresponding to intersections of  $\LLE \langle X_t 
\rangle$ with $\{Y_u\}$.

To bound $\E Z$,  
let $\widehat T_1$ and $\widehat T_2$ be independent
geometric random variables with mean 
$(1-\lambda)^{-1}=(\eps \vol)^{1/2}$,
and let $U$ denote the set of times $0\leq t \leq T_X$ that are locally
retained.
As $\langle X_t \rangle$ is symmetric, we can extend it to a doubly 
infinite chain $\langle X_t \rangle_{t=-\infty}^\infty$.
Since $Y_0$ is uniformly distributed,
\begin{align*}
\E J_{ij} & = \frac{1}{\vol} \P[i \in U] \\
& \geq \frac{1}{\vol}\left(\P[\langle X_t \rangle_{t=-\widehat 
{T_1}}^0 \cap 
\langle X_t \rangle_{t=1}^{\widehat {T_2}} = \emptyset]
- 2 \P[\widehat{T_1} \leq s] \right).
\end{align*}
Using Lemma \ref{cutpointlem1} and observing
that $\cc \geq 1$ yields
$$\E Z \geq \frac{T_X T_Y}{ \cc \vol} 
\left[\frac{1}{1+2\eps}-2(\sqrt{\eps})\right].$$

To bound the second moment, 
note that the number of intersections of the local 
loop-erasure starting at $x$ with the path 
from $y$ is bounded above by the number of intersections of the 
original walk, so it suffices to bound the second moment of the number 
of such intersections.
To do so, let $I_{ij}$ be an indicator random variable for the event 
$\{X_i=Y_j\}$ and fix a basepoint $o\in G_n$.
Because the starting points of our walks are uniform,
\begin{align*}
\E Z^2 & \leq \sum_{i=0}^{T_X} \sum_{j=0}^{T_Y} \sum_{k=i}^{T_X} 
\sum_{\ell=j}^{T_Y} 4\E I_{ij} I_{k\ell} - 3 \sum_{i=0}^{T_X} 
\sum_{j=0}^{T_Y} \E I_{ij} \\
&= \sum_{i=0}^{T_X} \sum_{j=0}^{T_Y} \sum_{k=i}^{T_X}
\sum_{\ell=j}^{T_Y} \frac{ 4 \P_o[X_{k-i+\ell-j}=o]}{\vol}
-\frac{3 T_X T_Y}{\vol} \\
& \leq \frac{4 T_X T_Y}{\vol} \sum_{k=0}^{T_X} \sum_{\ell=0}^{T_Y} 
\P_o[X_{k+\ell}=o] -\frac{3 T_X T_Y}{\vol} \\
& \leq [4\cc-3] \frac{T_X T_Y}{\vol}
\end{align*}
where $\cc$ is as in (\ref{localstrongtran}).
Using these quantities to lower bound $(\E Z)^2/ \E (Z^2)$ gives the 
desired result.
\end{proof}

Let $\Cap(S)$ be as in Definition \ref{capdef}.
Taking $T_X=q-3s$ and $T_Y=r$ gives:
\begin{corollary}
\label{expcapcor}
Suppose that $\{G_n\}$ is a sequence of graphs satisfying the assumptions 
of Theorem \ref{mainthm}.
If $q,r \leq \vol^{1/2-\delta/4}$, and $s=o(q)$, then
for simple random walk $\{X_t\}$ on $G_n$,
\begin{equation}
\label{expcapeqn}
\E\left[ \Cap(\LLE \langle X_t \rangle_{t=2s}^{q-s})  \right] 
\geq \frac{qr(1+o(1))}{\vol \cc^2(4\cc-3)}.
\end{equation}
\end{corollary}

\section{Tight concentration of longer segments}
\label{largedevsection}

We now combine the results of the previous section with large deviation 
estimates to show that $| \LLE(A_i)|$ and $\Cap [\LLE(A_i)]$ are 
tightly 
concentrated about their means.  We will do this by viewing $A_i$ as a 
union of smaller pieces and using the following 
large deviation bound:

\begin{lemma}[Hoeffding]
\label{largedevlem}
Suppose $\{Z_i\}$ are a family of independent random variables 
such that $0\leq Z_i \leq b$.  Then
\begin{equation}
\label{largedeveqn}
\P\left[\sum_{i=1}^{n} Z_i-\E Z_i > nt \right] \leq \exp
\left[-2n\left(\frac{t}{b}\right)^2 \right]. 
\end{equation}
\end{lemma}

This lemma is from
\cite{Hoeffding63}, Theorem 1, equation (2.3), and is proved using
standard arguments.

Our applications of this large deviation bound
include showing that, with high probability,
 the capacity and length 
of $\LLE(A_i)$ are close to their means for all $A_i$,
as well as the fact that for different $i$, the sequences
$\LLE(A_i)$ are
relatively far apart in the graph.  
Much of the difficulty in doing this for 
$\Cap[\LLE(A_i)]$ involves the possibility that a run of length $r$ of 
a walk might hit more than one subsegment of $\LLE(A_i)$.
These multiple intersections are important because they mean
that capacity is subadditive.

\begin{definition}
\label{closenessdef}
Let $\{Y_u\}$ be simple random walk on a graph $G$, 
and let $T_S=\min\{u\geq 0: Y_u \in S\}$ denote the hitting
time for $S$.  For sets $U$ and $V$,
the closeness of $U$ and $V$
is given by
$$
\close{U}{V}=\P_\pi[T_U<r, T_V < r].
$$
\end{definition}

Note that $\close{U}{V}$ will be small
 if $U$ and $V$ intersect in a 
single point yet are otherwise very far away.
On the other hand,  
if $U$ and $V$ are disjoint, but $V$ is a translation
by a small fixed distance, then $\close{U}{V}$ will be large.
Closeness is primarily a measure of whether or not typical points
of $U$ are near $V$, and vice versa, and is maximized when the two sets 
coincide.

In the case when $V$ is a segment of a random walk, 
$\close{U}{V}$ is a random variable whose mean is bounded 
by the following lemma:

\begin{lemma}
\label{closenesslem}
Let $\{G_n\}$ be a sequence of graphs satisfying the assumptions
of Theorem \ref{mainthm}, and let $r=r(n)$ be a sequence of positive 
integers.
Let $\{X_t\}$ be a random walk on $G_n$ started with $X_0=x$, 
let $T\geq 0$ be a random time that is independent of $\{X_t\}$,
and let $V=\{X_t\}_{t=\mix}^{T}$.  
Then for any set $U$ and starting position 
$x$,
\begin{equation}
\label{expcloseness}
\E \close{U}{V} \leq 4 \frac{r \Cap[U] \E T}{\vol}.
\end{equation}
\end{lemma}

\begin{proof}
Let $\{Y_u\}$ be a simple random walk on $G_n$ started in the stationary 
distribution $\pi$.
By the strong Markov property, for any fixed set $W$, 
\begin{multline*}
\close{U}{W} \leq \P_\pi[T_U< r] \P_{Y_{T_U}}[T_W < r
\mid T_U<r] \\
+ \P_\pi[T_W<r] \P_{Y_{T_W}}[T_U<r\mid T_W<r].
\end{multline*}
To apply this to our random $V$, let $\E_x$ denote
expectation given $X_0=x$.
\begin{align*}
\E_x & \left[\P_\pi[T_V<r] \P_{Y_{T_V}}[T_U<r \mid T_V<r]
\right]\\
 & = \E_x \left[\sum_y \P_\pi[T_V<r, Y_{T_V}=y] \P_y[T_U<r]\right] \\
& = \sum_y \P_y[T_U<r] \E_x\left[\P_\pi[T_V<r, Y_{T_V}=y]\right] \\
& \leq \sum_y \P_y[T_U<r]  \sum_{t=\mix}^{\infty} \sum_{u=0}^{r-1} 
\E_x \P_\pi[Y_u=X_t=y, t\leq T]\\
& < \sum_y \P_y[T_U<r]
  \frac{2r \E T}{\vol^2}\leq \frac{2r \E T \Cap[U]}{\vol}.
\end{align*}
Likewise, $\E_{x} \left[\P_{Y_{T_U}} [T_V < r 
\mid T_U<r] \right]$
is bounded by
\begin{align*}
& \sum_y \P_\pi[Y_{T_U}=y\mid T_U<r]
\sum_{u=0}^{r-1} \sum_{t=\mix}^{\infty} \P[Y_u=X_t, t\leq T
\mid Y_0=y, X_0=x] \\
\leq & \sum_y 
\P_\pi[Y_{T_U}=y\mid T_U <r] \frac{2 r \E T}{\vol}=
\frac{2 r \E T}{\vol}, 
\end{align*}
where the second inequality uses the fact that as
$t \geq \mix$, 
$$\P[Y_u=X_t \mid Y_0=y, X_0=x] \leq 2/\vol.$$
As $\P_\pi[T_U<r]=\Cap[U]$, this completes the proof.
\end{proof}

Let $\alpha=\alpha(n)$ and $\gamma=\gamma(n)$ be as in 
(\ref{alphagammadef}), implying that
$$
\E | \LLE\{X_t\}_{t=2s}^{r-s}|= \gamma r, \qquad
\E \Cap[ \LLE \{X_t\}_{t=2s}^{r-s} ] = \frac{\alpha r^2}
{\vol}.
$$

\begin{lemma}
\label{concentrationlem}
Let $\{G_n\}$ be a sequence of graphs satisfying the assumptions
of Theorem \ref{mainthm}
and let $\{X_t\}$ be simple random walk on $G_n$.  
Suppose that $r=r(n)$ and $s=s(n)$ are sequences such that
$$4 s^{1/2} \vol^{1/4} \log \vol\leq r \leq \vol^{1/2-\delta/4}$$
Taking $A_i$ as in (\ref{Aidef}), then 
whenever $\vol \geq 2^8$,
\begin{equation} 
\label{lengtheqn}
\P \left[ \Bigl|  |\LLE(A_i)| - \gamma r \Bigr| > 2 r\left(
\frac{s}{r}\right)^{1/6} 
\right]
\leq \exp\left[-2\left(\frac{r}{s}\right)^{1/6} \right],
\end{equation}
\begin{equation}
\label{capeqn}
\P\left[ \left|\Cap(\LLE(A_i))-\frac{\alpha r^2} {\vol} 
\right| > \left(\frac{r}{\sqrtvol}\right)^{9/4} 
\right]\leq 9 \left(\frac{r^{7/4}}{\vol^{7/8}}\right),
\end{equation}
Moreover,
$\alpha$ and $\gamma$ satisfy the bounds
\begin{equation}
\label{alphagammaeqn}
\gamma \geq \frac{1+o(1)}{\cc}, \qquad \alpha \geq 
\frac{1+o(1)}{\cc^2(4\cc-3)}.
\end{equation}
\end{lemma}

The situation that we are most interested in is when
$s$ and $r$ are as in Table \ref{deftable}, in which case
assumption (\ref{fastmixing}) shows that
Lemma \ref{concentrationlem} applies.
More generally, note that the assumptions on $s$ and $r$ imply that 
as $\vol \rightarrow \infty$, 
$s=o(r)$ and also $s \vol^{-1/2}=o(r^2 \vol^{-1})$.

\begin{proof}
Let $q=\lfloor (rs)^{1/2} \rfloor$.
To prove all three parts of this lemma, we will break the interval 
$A_i$ down into $\lfloor r/q\rfloor$ smaller pieces $B_{i,k}$ 
whose length is asymptotically $q$, 
and use Lemma \ref{largedevlem} and the results of 
Section \ref{lemmasection}.  Let $B_{i,k}=\Z \cap [ir+kq+2s, 
ir+(k+1)q-s]$, 
and denote
$$\LLE(B_{i,k})=\langle X_t \rangle_{t\in B_{i,k}\cap U},$$
where $U$ is the set of times that are locally retained.

We begin with (\ref{lengtheqn}).  Breaking $\LLE(A_i)$ 
into pieces and summing gives
\begin{align*}
\sum_{k=0}^{\lfloor r/q \rfloor  -1} | \LLE (B_{i,k})| & \leq 
| \LLE(A_i)| \\
& \leq \sum_{k=0}^{\lfloor r/q \rfloor}  | \LLE (B_{i,k} )|
 + 3s\left(\frac{r}{q}+1\right). 
\end{align*}
By our choice of $q$, $rs/q \sim r(s/r)^{1/2}$, 
which is of a lower 
order than $r(s/r)^{1/6}$.  
The restrictions on $r$ and $\vol$ insure that
$3s(r+q)/q \leq r(s/r)^{1/6}$.
The spacing between the $B_{i,k}$ is 
such that $|\LLE(B_{i,k})|$ are i.i.d, so
applying Lemma \ref{largedevlem}
with 
$Z_k=|\LLE(B_{i,k})|$, $b=q$, and $t=q^{4/3} r^{-1/3}$
shows that 
\begin{align*}
\P\left[\sum_{k=0}^{\lfloor r/q\rfloor -1} (Z_k-\E Z_k) > r(s/r)^{1/6} 
\right] & \leq \exp\left[-2\frac{r}{q} 
\left(\frac{q^{1/3}}{r^{1/3}}\right)^2 \right] \\
& \leq \exp \left[-2 \left(\frac{r}{s}\right)^{1/6} \right].
\end{align*}

The argument for (\ref{capeqn}) is similar, but the naive upper and lower 
bounds are farther apart.
The capacity is bounded above by the sum of the capacities 
of the pieces, plus a little extra since the $B_{i,k}$ are spaced 
$3s$ steps apart in time.  Summing the capacity 
of the pieces overcounts by the probability of hitting at least two 
pieces.  For a lower bound, we subtract the probability
of having double hits.
 \begin{multline*}
\sum_{k=0}^{\lfloor r/q \rfloor -1} \Cap[\LLE(B_{i,k})]
  - \sum_{k,j, k\neq j} \close{B_{i,k}}{B_{i,j}}  \leq \Cap[\LLE(A_i)] 
\\
\leq \sum_{k=0}^{\lfloor r/q \rfloor} \Cap[\LLE(B_{i,k})] + 
3\left(\frac{r}{q}+1\right) \frac{sr}{\vol}.
\end{multline*}
Applying
Lemma \ref{largedevlem}
with $Z_k=\Cap[\LLE(B_{i,k})]$, $b=qr/\vol$,
and $t=\left(qr^{5/4}\vol^{-9/8}\right)/2$
gives
$$
\P\left[\sum_k Z_k-\E Z_k 
> \frac{1}{2} \left(\frac{r}{\sqrtvol}\right)^{9/4} \right]\leq 
\exp\left[-\frac{1}{2} 
\frac{r^{3/2}}{s \vol^{1/4}} \right],
$$
which is at most $(1/2) \vol^{-2}\leq r^{7/4}\vol^{-7/8}$
by the assumptions on $r$.
By Lemma \ref{closenesslem}, 
$$\sum_{k,j} \E \close{B_{i,k}}{B_{i,j}}
\leq 4 \left(\frac{r}{q}+1\right)^2 \frac{q^2 r^2}{\vol^2}.$$  
By Markov's inequality,
$$
\P\left[\sum_{k,j, k\neq j} \close{B_{i,k}}{B_{i,j}} \geq \frac 
{r^{9/4}} {2\vol^{9/8}} 
\right] \leq \frac{17}{2}\left(r^{7/4} \vol^{-7/8} \right).
$$
This bound proves (\ref{capeqn}).

Finally, the lower bounds of (\ref{alphagammaeqn})
are those of (\ref{explengtheqn}) and (\ref{expcapeqn}).
\end{proof}

\section{Completion of the coupling argument}
\label{couplingsection}

The aim of this section is to prove Lemma \ref{couplinglemma}.
First, we will show that $\LE \langle X_t \rangle_{t=0}^T$ can be 
decomposed into runs of length $r$ without losing too much information, 
then we will prove two lemmas that allow us to couple random walk on 
$\extG{n}{L}$ with $\extK{m}{\widetilde L}$ for suitable $m$ and 
$\widetilde L$, and then we will end by proving Lemma 
\ref{couplinglemma}.

Given $T, r$ and $s$,
a time index $i \leq \ell=\lceil T/r\rceil$ is called {\it good} by time 
$T$ if 
$$
\{X_t\}_{t\in A_i} \cap \left[\{X_t\}_{t\leq ir} \cup \{X_t\}_{t\in
[(i+1)r, T]}\right] = \emptyset,
$$
where $A_i$ is as in (\ref{Aidef}).  Intuitively, this
means that there are no loops longer than 
length $s$ with one endpoint inside $A_i$ and one outside, 
although there are some slight differences at times near the endpoints of 
$A_i$.  Similarly, a time index 
$i\leq \ell$ is called a {\it single intersection} 
at time $T$ if there exists $j\leq \ell$ such that $\{X_t\}_{t\in A_i} 
\cap 
\{X_t\}_{t\in 
A_j}\neq \emptyset$, but
$
\{X_t\}_{t\in A_i} \cap \{X_t\}_{t\in I} = 
\emptyset,
$
where $I=\{[0,T] \setminus([(i-1)r, ir] \cup [(j-1)r,jr] ) \}\cap \Z$.

Finally, a time index is called {\it bad} if it is neither
a single intersection nor good.
Let $\bad{T}$ denote the collection of time indices $i\leq \ell$ 
that are bad at time $T$, let $\single{T}$ be those that are 
single intersections, and $\good{T}$ be those that are good.

\begin{definition}
\label{locallydecompdef}
Given $T, r$ and $s$, let $\ell=\lceil T/r \rceil$.
A run $\langle X_t \rangle_{t=0}^T$ is called
{\it locally decomposable} if:
\begin{enumerate}
\item
$
\bad{T}=\emptyset$
\item
\label{tightlength}
for all $i\leq \ell$, 
$\Bigl| | \LLE(A_i)|-\gamma r \Bigr| \leq 2 r 
\left(\frac{s}{r}\right)^{1/6}$
\item
\label{tightcap}
for all $i\leq \ell$, 
$| \Cap[\LLE(A_i)]-\alpha r^2 \vol^{-1}| \leq r^{9/4} \vol^{-9/8}$
\item
\label{nointermedloops}
There is no pair $(t,u)\in \Z^2$ with $t,u\in [0,T]$ and $|t-u|\in 
[\mix,r]$
such that $X_t=X_u$.
\item
\label{nohitgaps}
There is no pair $(t,u)\in \Z^2$ with $t,u \in [0,T]$ and $|t-u|\geq 
s$ such that $t\notin \cup_i A_i$ and $X_t=X_u$.
\item
For all $t\in [0,T-s]$, the interval $[t,t+s]$ contains a local
cutpoint.
\end{enumerate}
\end{definition}

\begin{lemma}
\label{locallydecomplemma}
Let $\{X_t\}$ be simple random walk on $G_n$ and $T$ a 
geometric random variable with mean $L \sqrtvol$ that is independent
of $\{X_t\}$.  Suppose that $r$ and $s$ satisfy the hypotheses of Lemma 
\ref{concentrationlem}, and also that $s>\mix (\log \vol)^2$.
Then the probability that $\langle X_t \rangle_{t=0}^T$
is locally decomposable is $1-O(L^2 r^{3/4}\vol^{-3/8})$.
Moreover, 
\begin{equation}
\label{expnotgood1}
\E[\bad{T}\cup \single{T}] =O(L^2).
\end{equation}
\end{lemma}

When $s,q$, and $r$ are as in Table \ref{deftable}, 
condition (\ref{fastmixing}) implies that Lemma \ref{locallydecomplemma} 
applies and also that $r^{3/4}\vol^{-3/8}=o(\vol^{-3\delta/16})$.

\begin{proof}
Counting the expected number of long loops involving
$\LLE(A_i)$ gives
$$\P[ i\in \bad{T}\cup \single{T} \mid 
T] \leq \frac{2 T r}{\vol}.$$
As $\ell \leq (T/r)+1$, this implies that
\begin{equation}
\label{expnotgood}
\E[|\bad{T} \cup \single{T}| \mid T] 
\leq \frac{2 (T^2+rT)}{\vol}.
\end{equation}
Using the restrictions on $r$ and
taking the expected value of (\ref{expnotgood}) yields
(\ref{expnotgood1}).

Likewise, the probability that an index $i$ is bad is bounded by the 
expected number of pairs $j$ and $k$ such 
that $\{X_t\}_{t \in A_i}$ intersects both $\{X_t\}_{(j-1)r\leq t
\leq jr}$ and 
$\{X_t\}_{(k-1)\leq t \leq k}$, 
plus the expected number of times such that $\{X_t\}_{t\in A_i}$ 
intersects the gaps between the $A_j$,
yielding
\begin{align*}
\P[i \in \bad{T} \mid T] & \leq 
\binom{\lceil T/r \rceil}{2} 
\left( \frac{2 r^2}{\vol}\right)^2
+ \left \lceil \frac{T}{r} \right\rceil
\frac{6sr}{\vol} \\ 
& \leq \frac{2 (T+r)^2 r^2}{\vol^2}+ \frac{6(T+r)s}{\vol},
\end{align*}
and as $\E T = L \sqrtvol$
and $s \vol^{-1/2}=o(r^2 \vol^{-1})$,
 we see that 
$\E |\bad{T}| = O\left(L^3 r \vol^{-1/2}\right).$  
In particular, $\P[\bad{T}=\emptyset]=1-O(L^3 r \vol^{-1/2})$.

Consider the
sequence of events $\step{i}$ given by:
$$
\step{i}=
\left\{
\begin{array}{c c}
\Bigl| | \LLE(A_j)|-\gamma r\Bigr| \leq 2r\left(\frac{s}{r}\right)^{1/6}, 
& j\leq i \\
\Bigl| \Cap(\LLE(A_j)) - \alpha r^2 \vol^{-1} \Bigr| \leq
r^{9/4}\vol^{-9/8}, & j \leq i \\
\end{array}
\right\}.
$$
Using Lemma \ref{concentrationlem} and the facts that $\ell 
=\lceil T/r \rceil$ and $\E T=L \sqrtvol$,
$$
\P[\step{\ell} \mid T] \geq 1- 10 \frac{T+r}{r}
\frac{r^{7/4}}{\vol^{7/8}},
$$
$$
\P[\step{\ell}]=1-O\left(L\frac{r^{3/4}}
{\vol^{3/8}}\right).
$$
Let $\Gamma_n$ denote the event that for 
all $t\in [0,T-s]$, the interval $[t,t+s]$ contains a local cutpoint.
By Corollary \ref{cutptcor}, there exists $C$ so that
 $\P[\Gamma_n]=1-L \sqrtvol \exp(-Cs/\mix)$, which is 
$1-O(Lr^{3/4}\vol^{-3/8})$ by the lower bound on $s$.

To bound the probability of Condition \ref{nointermedloops} of 
Definition \ref{locallydecompdef} holding,
conditioned on $T$,
the expected number of loops of length in the interval $[\mix, r]$
is
$$
\sum_{i=0}^{T} \sum_{j=i+\mix}^{i+r} \P[X_i=X_j]\leq 
\frac{2T r}{\vol}.
$$
As $\E T= L \sqrtvol$,
Condition \ref{nointermedloops} holds with probability 
$1- O(Lr/\sqrtvol)$.

For Condition \ref{nohitgaps}, conditioning on $T$ and
counting 
the expected number of pairs $t$ and $u$ that we wish to
avoid bounds the probability of failing by
$$
3s \left(\frac{T}{r}+1\right) \frac{2T}{\vol},
$$
so Condition \ref{nohitgaps} holds with probability
at least $1-O(sL^2/r)$.  Combining all six of these bounds proves the 
claim.
\end{proof}

For a finite tree $\T$, call a vertex $v\in \T$ a {\it leaf} if the 
degree of $v$ is 1.

\begin{lemma}
\label{uniformtreelemma}
Suppose $\T$ and $\widehat \T$ are trees with $k$ leaves $\{x_1, \dots, 
x_k\}$ 
and $\{y_1, \dots, y_k\}$ respectively.  Let $\xi$ and $\eta$ be uniformly 
chosen from $\T$ and $\widehat \T$ respectively.
There exist constants $C_1, C_2$ depending on $k$ such that if
for some $\eps>0$ and $R\geq 1$,
 $|d(x_i,x_j)-R d(y_i,y_j)|\leq \eps d(x_1, x_2)$ for all $i,j\leq k$, 
then there exists 
 a coupling of $\xi$ and $\eta$ such that
\begin{equation}
\label{uniformtreeeqn}
\P\Bigl[\forall \, i\leq k, |d(\xi,x_i)-R d(\eta, y_i)| \leq C_1 [\eps
d(x_1,x_2)+R]
\Bigr] \geq 1-C_2 \eps.
\end{equation}
\end{lemma}

\begin{proof}
In the case when $k=2$, take $\mu$ to be a uniform random variable on 
$[0,1]$, take $\xi$ so that $d(\xi, x_1)=\lfloor \mu[d(x_1,x_2)+1] 
\rfloor$, and take $\eta$ so that $d(\eta, y_1)=\lfloor 
\mu[d(y_1,y_2)+1]\rfloor$.  We then have $|R d(\eta, y_1)-d(\xi, x_1)| 
\leq \eps d(x_1, x_2) + R +1$.  As $R\geq 1$,
this shows that we can take $C_1=2$ and 
$C_2=0$.

Proceeding by induction, assume that we know the result for $k=n$.  Given 
a tree with $n+1$ ends, the vertex $x_{n+1}$ is connected to the subtree 
spanned by $\{x_1, \dots, x_n\}$ by a segment of length 
$$A=\min_{i,j\leq n} [d(x_i,x_{n+1})+d(x_j, x_{n+1})-d(x_i, x_j)].$$
Likewise, let $\widehat A$ denote the length of the segment connecting 
$y_{n+1}$ to the tree spanned by $\{y_1, \dots, y_n\}$.  
By repeated application of the triangle inequality, our assumption on 
the distances shows that $|A -R \widehat A| \leq 6 \eps d(x_1, x_2)$ and
also that
$$
\Bigl|\P[d(\xi, x_{n+1}) <A]-\P[d(\eta, y_{n+1})<\widehat A] \Bigr|
\leq C \eps
$$ 
for some constant $C=C(n)$.
If both $\xi$ and $\eta$ are chosen on the segment connecting $x_{n+1}$ 
and $y_{n+1}$ to the rest of the tree, then we can couple them so that
$|d(\xi, x_{n+1})-R d(\eta, y_{n+1})| \leq 2R+6 \eps d(x_1, x_2)$, and so 
by 
the triangle inequality, for all $i\leq n$,
$|d(\xi, x_{i})-R d(\eta, y_{i})| \leq 2R+8 \eps d(x_1, x_2)$.
If both are in the rest of the tree, then the induction hypothesis shows 
that with probability $1-C_2(n) \eps$, we can couple them so that
$|d(\xi, x_{i})-R d(\eta, y_{i})| \leq C_1(n)[\eps d(x_1, x_2)+R]$
for all $i\leq n$, and thus by the triangle inequality,
$|d(\xi, x_{n+1})-R d(\eta, y_{n+1})| \leq (C_1(n)+6) \eps d(x_1, 
x_2)+C_1(n) R$.  We have thus 
proved the lemma with 
$C_1(n+1)=C_1(n)+6$, $C_1(1)=2$, 
and $C_2(n+1)=C_2(n) + C(n)$, $C_2(1)=0$.
\end{proof}

We now begin constructing the coupling.
Let $\extG{n}{L}$ and $\extK{m}{\widetilde L}$ be
as in Lemma \ref{couplinglemma}, with
$m$ and $\widetilde L$ given by
\begin{equation}
\label{mdef}
m= \left \lfloor \frac{\vol}{\alpha r^2} \right \rfloor
\qquad \mbox{\rm and}
\end{equation}
\begin{equation}
\label{tildeLdef}
\widetilde L =\frac{1}{\sqrt{m}}\left[1-\left(1-\frac{1}
{L \sqrtvol}\right)^r \right]^{-1}.
\end{equation}
For fixed $L$, note that
$\widetilde L \alpha^{-1/2}\rightarrow L$ as
$n\rightarrow \infty$ since $r=o(\sqrtvol)$.

\begin{lemma}
\label{inductionsteplemma}
Suppose that $\{G_n\}$ is a sequence of graphs satisfying the assumptions 
of Theorem \ref{mainthm}, and $r=r(n)$ and $s=s(n)$ are sequences
of constants satisfying the assumptions of Lemma \ref{locallydecomplemma}.
For fixed $L$, 
let $\{X_t\}$ be weighted random walk on $\extG{n}{L}$
with $X_0$ chosen uniformly from $G_n$.
Taking $m$ and $\widetilde L$
as in (\ref{mdef}) and (\ref{tildeLdef}), 
let $\{Y_u\}$ be weighted random walk on $\extK{m}{\widetilde L}$
with $Y_0 \in K_m$, and 
let $\Kindex{i}$ and $\Gindex{i}$ be as in (\ref{Kindexdef}) and 
(\ref{Gindexdef}).
Suppose that $\T_k\subset \extG{n}{L}$ and 
$\widetilde \T_k \subset \extK{m}{\widetilde \T_k}$ are such that $\rho 
\in \T_k$ and $\tilde \rho \in \widetilde \T_k$.
Then there exists a coupling of $\{X_t\}$ and $\{Y_u\}$
such that
\begin{multline*}
\P\left[\widetilde T = \left \lceil\frac{T}{r}\right 
\rceil,
\Gindex{j}= \Kindex{j} \, \forall \, j\leq 
\widetilde T
\right]=1-O\left(\frac{L^2 r^{1/4}}{\vol^{1/8}}\right)
-O\left(\frac{Ls |\T_k|}{r \sqrtvol}\right) \\
- \frac{L \sqrtvol}{r} \Bigl| \Cap[\T_k \cap G_n]-\frac{|\widetilde \T_k 
\cap K_m|}{m}
\Bigr|
,\end{multline*}
where $T$ is the hitting time of $\T_k$ and $\widetilde T$
is the hitting time of $\widetilde \T_k$.
\end{lemma}

Note that for $s,q$ and $r$ as in Table \ref{deftable},
when (\ref{fastmixing}) and (\ref{couplecap}) hold, the hypotheses are 
met and also that the bound obtained is $1-O(\vol^{-\delta/32})$.

\begin{proof}
To establish the coupling, we will use the fact that while the 
distribution of $\Gindex{j}$ depends on more than just $\Gindex{j-1}$, the 
fact that $Y_u$ is uniformly chosen from $K_m$ means that the
distribution of $\Kindex{j}$ depends only on the value of $\Kindex{j-1}$.  
To this end, we will run the process $\{X_t\}$ on $\extG{n}{L}$, use 
that to 
evaluate the indicator random variables $\{\Iij\}_{i\in \Gindex{j-1}}$ 
(and hence compute 
$\Gindex{j}$), couple the indicators $\{\Iij\}_{i\in \Gindex{j-1}}$ 
with $\{\tildeIij\}_{i\in \Gindex{j-1}}$, and use the 
fact that we can construct $Y_j$ given the values of the
indicators $\{\tildeIij\}_{i\in \Kindex{j-1}}$.

We need to show two things: first, that we can couple
$\lceil T/r \rceil$ and $\widetilde T$, 
and second, that we can couple the walks until these hitting times.  

To couple the hitting times, 
let $T_\rho$ denote the first time that 
$X_t=\rho$, and note that $T_\rho \geq T$.  Until time 
$T_\rho$, the path of $\{X_t\}$ is a simple random walk on $G_n$ that is 
independent of $T_\rho$, so
a first moment estimate gives
$$\P[T \notin \cup A_i \mid T_\rho]
\leq \frac{|\T_k| 3s}{\vol} \left\lceil \frac{T_\rho}{r} \right \rceil.
$$
As $\E T_\rho =L \sqrtvol$, this yields
$
\P[T \notin \cup A_i]= O \left(L s |\T_k|r^{-1} \vol^{-1/2}\right).
$

Likewise, let $\widetilde T_\rho$ denote the
hitting time of $\widetilde \rho$.
By the definition of $\widetilde L$,
\begin{align*}
\P[\widetilde T=j & \mid \widetilde T>j-1]  =
\frac{|\widetilde \T_k \cap K_m|}{m}\left(1-\frac{1}{\widetilde L 
\sqrt{m}}\right) + \frac{1}{\widetilde L \sqrt{m}}\\
& = \frac{|\widetilde \T_k \cap K_m|}{m}\left(1-\frac{1}{L
\sqrtvol}\right)^r +1-\left(1-\frac{1}{L \sqrtvol}\right)^r
\end{align*} 
To estimate what happens on $\extG{n}{L}$, we will
use the fact that $s$ is much larger than $\mix$.
In particular, as the $L^\infty$ distance
between the distribution of $\{X_t\}$ and the
stationary distribution $\pi$, given by
$\sup_x \Bigl| p^t(x) \vol^{-1}-1\Bigr|$, is submultiplicative
(see e.g., \cite{SC04} Proposition 2.2), we
have $\P_x[X_s=y]=[1+O(2^{-s/\mix})]/\vol$, which is 
$[1+o(\vol^{-1})]/\vol$ by the lower bound on $s$.
Using this to correct the fact that our
walk is not quite uniform in distribution gives
\begin{align*}
& \Capr{r-3s}(\T_k \cap G_n)[1+O(2^{-s/\mix})] + 1-\left(1-\frac{1}{L 
\sqrtvol}\right)^{r-3s}  \\
\geq & \P[\langle X_t \rangle_{t\in A_j} \cap 
\T_k\neq \emptyset \mid T> (j-1)r-s]\\
\geq & \left(1-\frac{1}{L \sqrtvol}\right)^{3s} \left[ 
\frac{\Capr{r-3s}(\T_k \cap G_n)}{1+O(2^{-s/\mix})}+ 
1-\left(1-\frac{1}{L
\sqrtvol}\right)^{r-3s}\right].
 \end{align*}
Combining the estimates on $\extG{n}{L}$ and $\extK{m}{\widetilde L}$, 
we obtain
\begin{multline*}
\Bigl| \P[\langle X_t \rangle_{t\in A_j} \cap
\T_k\neq \emptyset \mid T> (j-1)r-s] - \P[\widetilde T=j \mid 
\widetilde T>j-1] \Bigr|\\
\leq O(2^{-s/\mix}) + O\left(\frac{s}{L \sqrtvol}\right) + \Bigl| 
\Capr{r-3s}(\T_k 
\cap G_n)- \frac{|\widetilde T_k \cap K_m|}{m} \Bigr|.
\end{multline*}

Note that 
$\Cap[\T_k]-\Capr{r-3s}[\T_k]=O(s|\T_k| \vol^{-1})$ 
by considering the expected number of intersections in the 
interval $[r-3s,r]$.
  Using the lower bound on $s$ to bound $2^{-s/\mix}$
and the fact that $T$ is bounded by a geometric random variable with
mean $L \sqrtvol$,
we have thus shown that the hitting times can be coupled with the claimed 
probability.

We now show that, with high probability, that
the two walks are coupled up until time
$T$.  Consider the event that
the first step at which the coupling breaks
is step $j$, which in turn is bounded by
\begin{equation}
\label{firstjprob}
\P[\Gindex{j-1}=\Kindex{j-1}, \Gindex{j}\neq \Kindex{j}].
\end{equation}

To bound (\ref{firstjprob}), we will bound the
total variation distance between the joint distributions of 
$\{\Iij\}_{i\in \Gindex{j-1}}$ 
and $\{\tildeIij\}_{i\in \Kindex{j-1}}$ conditioned
on the event $\Kindex{j-1}=\Gindex{j-1}$
as well as the past of the walk $\{X_t\}_{t<r(j-1)}$.  
To do so, 
consider three types of cases for the values
of $\{\Iij\}_{i \in \Gindex{j-1}}$.

Let Case 1 be the case when no long loop is formed
at step $j$, meaning that we need to 
bound
\begin{equation}
\label{case1bound}
\left| \P[ \Iij=0 \, \forall \, i \in \Gindex{j-1}]
-  \P[ \tildeIij=0 \, \forall \, i \in \Gindex{j-1}]
\right|.
\end{equation}
Conditioned on $\Gindex{j-1}=\Kindex{j-1}$,
$$\P[ \tildeIij=0 \, \forall \, i \in 
\Kindex{j-1}]=1-\frac{|\Gindex{j-1}|}{m}.$$
Again using a factor of $1+O(2^{-s/\mix})=1+o(\vol^{-1})$ to
correct for the fact that our
walk is not quite uniform in distribution,
\begin{align*}
1-& \left(\sum_{i\in \Gindex{j-1}} \Capr{r-3s}[\LLE(A_i)] \right)
\leq \P[ \Iij=0 \, \forall \, i \in \Gindex{j-1}]
[1+o(\vol^{-1})] \\
\leq & 1-\left(\sum_{i\in \Gindex{j-1}} \Capr{r-3s}[\LLE(A_i)] \right)
+ \sum_{i,k\in \Gindex{j-1}, i\neq k} \lclose{i}{j}.
\end{align*}
Combining these and using the triangle inequality, 
we see that (\ref{case1bound}) is,
up to a factor of $1+o(\vol^{-1})$, bounded by
\begin{equation}
\label{case1tvbound}
\sum_{i,k\in \Gindex{j-1}, i\neq k} \lclose{i}{j}
+ \sum_{i\in \Gindex{j-1}} \left| 
\Capr{r-3s}[\LLE(A_i)]- \frac{1}{m}\right|.
\end{equation}
Let Case 2 be the case when there is a long 
loop formed involving exactly one $A_i$, meaning that we need to bound
\begin{equation}
\label{case2bound}
\left| \P[\tildeIij=1, \tildeIkj=0 \, k\in \Kindex{j-1}
\setminus\{i\}]
- \P[\Iij=1, \Ikj=0\, k\in \Gindex{j-1}\setminus\{i\}]
\right|.
\end{equation}
But
\begin{align*}
\Capr{r-3s}[\LLE(A_i)] & \geq  \P[\Iij=1, \Ikj=0\, k\in 
\Gindex{j-1}\setminus\{i\}] 
(1+o(\vol^{-1}))\\ 
 \geq & \Capr{r-3s}[\LLE(A_i)] - \sum_{k\in \Gindex{j-1}\setminus\{i\}} 
\lclose{i}{k}
\end{align*}
and
$ 
\P[\tildeIij=1, \tildeIkj=0 \, k\in \Kindex{j-1}
\setminus\{i\}]=\frac{1}{m},
$
so (\ref{case2bound}) is,
up to a factor of $1+o(\vol^{-1})$, bounded by
\begin{equation}
\label{case2tvbound}
\left| \Capr{r-3s}[\LLE(A_i)] - \frac{1}{m} \right| + \sum_{k\in 
\Gindex{j-1}\setminus \{i\}} 
\lclose{i}{k}.
\end{equation}
Let Case 3 be the case that two or more
of the $\{\Iij\}_{i\in \Gindex{j-1}}$ are 1.
The probability of this is 0 on $K_m$ and on $G_n$ is at most
\begin{equation}
\label{case3tvbound}
\sum_{i,k \in \Gindex{j-1}, i\neq k} \lclose{i}{k}.
\end{equation}
Summing (\ref{case2tvbound}) over all $i$ in 
$\Gindex{j-1}$ and adding to 
(\ref{case1tvbound}) and (\ref{case3tvbound}) 
gives an upper bound for the 
total variation distance between the joint distributions of 
$\{\Iij\}_{i\in \Gindex{j-1}}$ and 
$\{\tildeIij\}_{i\in \Kindex{j-1}}$ 
conditioned on $\Gindex{j-1}=\Kindex{j-1}$ 
and the path $\{X_t\}_{t<(j-1)r}$ of, up to a factor of $1+o(1)$,
$$
3 \sum_{i,k \in \Gindex{j-1}, i\neq k} \lclose{i}{k}
+ 2 \sum_{i\in \Gindex{j-1}} \left| \Capr{r-3s}[\LLE(A_i)] 
- \frac{1}{m} \right|.
$$
Taking expectation to remove the conditioning,
$\Capr{r-3s}[\LLE(A_i)]\leq r^2/\vol$ 
by equation (\ref{easycapbound}). Using
Lemma \ref{closenesslem},
  $\E \lclose{i}{k}\leq 4r^4/\vol^2$.  Since $T_\rho$
  is independent of the path of $\{X_t\}$ for $t\leq T_\rho$, 
\begin{align*}
\E \sum_{j=1}^{\lceil T_\rho/r\rceil}  \sum_{i,k \in \Gindex{j-1}, i\neq 
k} 
\lclose{i}{k} & \leq \frac{4 r^4}{\vol^2} \frac{ \E(T_\rho+r)^3}{r^3}\\
& = O\left(\frac{r L^3}{\sqrtvol}\right).
\end{align*}
As $\Bigl|\alpha (r-3s)^2 \vol^{-1}-m^{-1}
\Bigr|=O(r^4 \vol^{-2})$,
Lemma \ref{concentrationlem} and the bound
from (\ref{easycapbound}) show that
\begin{align*}
\E \left| \Capr{r-3s}[\LLE(A_i)]-\frac{1}{m} \right| & \leq
\left(\frac{r}{\sqrtvol}\right)^{9/4} + 9\frac{r^2}{\vol}
\frac{r^{7/4}}{\vol^{7/8}} \\ & = 
O\left(\frac{r}{\sqrtvol}\right)^{9/4}.
\end{align*}
This leads to the bound
\begin{align*}
\E \sum_{j=1}^{\lceil T_\rho/r\rceil}  \sum_{i \in \Gindex{j-1}}
\left| \Capr{r-3s}[\LLE(A_i)] - \frac{1}{m} \right| &
\leq \frac{\E T_\rho^2}{r^2} O\left(\frac{r}{\sqrtvol}\right)^{9/4}
\\ & = O\left(\frac{L^2 r^{1/4}}{\vol^{1/8}}\right),
\end{align*}
which completes the proof.
\end{proof}

We now use Lemmas 
\ref{locallydecomplemma} -- \ref{inductionsteplemma}
to prove the induction
step of Lemma \ref{couplinglemma}.

\begin{proof}[Proof of Lemma \ref{couplinglemma}]
We proceed by induction.
Let $\{X_t\}$ be a weighted random walk on $\extG{n}{L}$
with $X_0=x_{k+1}$, and let 
$T$ denote the hitting time of $\T_k$.
Take $q,r$ and $s$ as in Table \ref{deftable},
and assume that Lemma \ref{couplinglemma} holds for $k$.
As before, take $\ell=\lceil T/r \rceil$.  If $i \in \good{T} \cap
\Gindex{\ell}$, 
and the path $\langle X_t \rangle_{t=0}^T$ is locally
decomposable,  then 
$\LLE(A_i) \subset  
\LE \langle X_t\rangle_{t=0}^T$ and, in particular,
\begin{equation}
\label{llelowerbound}
\sum_{i\in \Gindex{\ell}
\cap \good{T}} |\LLE(A_i)| \leq |\LE\langle X_t \rangle_{t=0}^T|.
\end{equation}
Moreover, when $\langle X_t \rangle_{t=0}^T$ is locally
decomposable, the only difference
between $\cup_{i\in \Gindex{\ell}} \LLE(A_i)$ and
$\LE \langle X_t \rangle_{t=0}^T$
are from the gaps of length $3s$ between the $A_i$, the 
runs that were erased during the single
intersections $\single{T}$, and the final run $\LLE \langle 
X_t\rangle_{t=(\ell-1)r}^T$.
 Thus, when $\langle X_t \rangle_{t=0}^T$ is locally
decomposable and $\Gindex{\ell}=\Gindex{\ell-1} \concat \{\ell\}$,
\begin{equation}
\label{lleupperbound}
 |\LE\langle X_t \rangle_{t=0}^T|\\
\leq \sum_{i\in \Gindex{\ell}}
|\LLE(A_i)| +
\left(\frac{3s T}{r}\right) + r(|\single{T}|+1).
\end{equation}
Combining (\ref{llelowerbound}) and (\ref{lleupperbound}),
$$
\Bigl|
 |\LE\langle X_t \rangle_{t=0}^T|-\gamma r |\Gindex{\ell}| 
\Bigr| \leq \left(\frac{3s T}{r}\right) + r(|\single{T}|+1)
+ \sum_{i\in \Gindex{\ell}} \Bigl| |\LLE(A_i)| - \gamma r \Bigr|.
$$

By the definitions
of Table \ref{deftable}, equation (\ref{expnotgood1}),
and condition (\ref{fastmixing}), 
$$\E\left[\frac{3s T}{r}\right] 
=O\left(L \mix^{1/2}\vol^{1/4}\right)=o\left(\vol^{1/2-\delta/2}\right)
\qquad \mbox{\rm{and}}$$ 
$$\E[r|\single{T}|]=O(L^2 r)=
o\left(\vol^{1/2-\delta/4}\right).$$  
When $\langle X_t \rangle_{t=0}^T$ is locally decomposable,
 $$\sum_{i\in \Gindex{\ell}}
 \Bigl| |\LLE(A_i)|-\gamma r \Bigr|
\leq 2\left[\frac{T+r}{r}\right]r\left(\frac{s}{r}\right)^{1/6}
=o\left(T \vol^{-\delta/12}
\right).$$
As $T$ is geometric with mean $L \sqrtvol$,
for large enough $n$, this is with probability $1-o(\vol^{-1})$
less than $\vol^{1/2-\delta/24}$.
Note that $$\P[\Gindex{\ell}= \Gindex{\ell-1} \concat \{\ell\} \mid T]
\geq Tr \vol^{-1},$$
and applying this along with Lemma \ref{locallydecomplemma} and Markov's 
inequality gives
\begin{equation}
\label{newLEeqn}
\P\left[ \left| \frac{ | \LE \langle X_t \rangle_{t=0}^T| - \gamma r 
|\Gindex{\ell}|}{\beta \sqrtvol}\right| > \vol^{-\delta/24} 
\right]=O(\vol^{-3\delta/16}).
\end{equation}
To compare with the complete graph,
$$
\left| \frac{|\LE\langle X_t\rangle_{t=0}^T|}
{\beta \sqrtvol} - \frac{|\Gindex{\ell}|}{\sqrt{m}}\right|
\leq
\left| \frac{|\LE\langle X_t\rangle_{t=0}^T| - \gamma r |\Gindex{\ell}|}
{\beta \sqrtvol}\right| 
+
|\Gindex{\ell}| \left|\frac{1}{\sqrt{m}}-\frac{r \alpha^{1/2}}
{\sqrtvol} \right|
$$
But $[m^{-1/2} - r\alpha^{1/2} \vol^{-1/2}]=O(r \vol^{-1})$,
and by Lemma \ref{inductionsteplemma}, with probability 
$1-O(\vol^{-\delta/32})$, 
$|\Gindex{\ell}|=|\Kindex{\ell}|$ is the length of the segment
added to $\widetilde \T_k$ to form $\widetilde \T_{k+1}$.
Let $\xi$ be chosen uniformly from $\T_k$ and according
to the coupling of (\ref{uniformhitting}).
Then
$$
d(x_i,x_{k+1})=d(x_i, \xi) + |\LE\langle X_t \rangle_{t=0}^T|
+ d(x_i, X_T)-d(x_i, \xi).
$$
Using the decomposition
\begin{multline*}
\left| \frac{d(x_i, x_{k+1})}{\beta \sqrtvol} - \frac{d(y_i, 
y_{k+1})}{\sqrt{m}} \right| \leq 
\left| \frac{|\LE\langle X_t\rangle_{t=0}^T|}
{\beta \sqrtvol} - \frac{|\Gindex{\ell}|}{\sqrt{m}}\right| \\
+ \left| \frac{d(x_i, \xi)}{\beta \sqrtvol} -\frac{d(y_i, Y_T)}{\sqrt{m}} 
\right| 
+ \frac{d(x_i, X_T)-d(x_i, \xi)}{\beta \sqrtvol},
\end{multline*}
we have already controlled the first term,
equation (\ref{uniformhitting}) implies that
$$\P[| d(x_i, X_T)-d(x_i, \xi)| > \vol^{1/2-\delta/32}]
\leq O(\vol^{-\delta/32}),$$ and
applying Lemma \ref{uniformtreelemma} to couple
$\xi$ and $Y_T$, we see that
there is a coupling such that (\ref{couplelength})
holds with probability $1-O(\vol^{-\delta/32})$.

To prove (\ref{couplecap}), 
we first bound
$\close{\T_k}{\LE\langle X_t \rangle_{t=0}^T}.$
By Lemma \ref{closenesslem},
\begin{equation}
\label{Tkcloseness}
\E [\close{\T_k}{ \langle X_t\rangle_{t=0}^{\killT}} \mid \T_k] 
\leq 4 \frac{L r \Cap{\T_k}}{\vol^{1/2}}.
\end{equation}
Markov's inequality gives
$$
\P\left[ \close{\T_k}{ \langle X_t\rangle_{t=0}^{\killT}} 
> 4 \frac{\Cap(\T_k) r^{1/2}}{\vol^{1/4}} \mid \T_k \right]
\leq \frac{L r^{1/2}}{\vol^{1/4}}.
$$
The fact that
 $\LE \langle X_t \rangle_{t=0}^T \subset \{X_t\}_{t=0}^{\killT}$
and monotonicity of closeness, along with the fact that
 $r^{1/2} \vol^{-1/4}=o(\vol^{-\delta/8})$, implies that
$$
\P[\close{\T_k}{\LE\langle X_t \rangle_{t=0}^T}\leq \vol^{-\delta/8}
\Cap(\T_k) \mid \T_k]=1-o(\vol^{-\delta/8}).
$$
Now estimating $\Cap[\LE\langle X_t \rangle_{t=0}^T]$,
when $\langle X_t \rangle_{t=0}^T$ is locally decomposable,
\begin{align*}
\sum_{i\in \Gindex{\ell}
\cap \good{T}} \Cap[\LLE(A_i)] - &
\sum_{i,j \in \Gindex{\ell}}
\lclose{i}{j} \\
& \leq \Cap[\LE\langle X_t 
\rangle_{t=0}^T]\\
& \leq \sum_{i\in \Gindex{\ell}}
[\Cap\LLE(A_i)] +
\left(\frac{3s T}{\vol}\right) + \frac{r^2(|\single{T}|+1)}{\vol}.
\end{align*}
In particular, we have
\begin{multline*}
\left| \Cap[\LE\langle X_t \rangle_{t=0}^T] 
- \frac{|\Gindex{\ell}|}{m} 
\right|
\leq \sum_{i \in \Gindex{\ell}} \left| \Cap[\LLE(A_i)]
- \frac{1}{m} \right| \\ +
\sum_{i,j\in \Gindex{\ell}} \lclose{i}{j}
+ \left(\frac{3s T}{\vol}\right)+\frac{r^2(|\single{T}|+1)}{\vol}.
\end{multline*}
The expected value of the sum of the last three of
these terms is bounded by
$$
\frac{4 r^2 \E (\killT)^2}{\vol^2}+
 \frac{3s L}{\sqrtvol}
+ \frac{2r^2 [\E (\killT)^2+1]}{\vol^2}
=O\left(\frac{r^2}{\vol}\right).
$$
By Markov's inequality, the sum of these three terms is 
therefore less than $\vol^{-\delta/16}$ with probability 
$1-o(\vol^{-\delta/32})$.  On the event $\step{\ell}$,
$$
\sum_{i \in \Gindex{\ell}} \left| \Cap[\LLE(A_i)]
- \frac{1}{m} \right| \leq \frac{\killT+r}{r}
\left[\frac{r^{9/4}}{\vol^{9/8}}+O \left(\frac{r^4}{\vol^2}
\right)\right],
$$
which is less than $\vol^{-\delta/20}(r/\sqrtvol)$
with probability $1-o(\vol^{-1})$
because $\killT$ is geometric with mean $L \sqrtvol$.
By Lemma \ref{locallydecomplemma}, the values
of Table \ref{deftable} and Condition (\ref{fastmixing}),
$\LE\langle X_t \rangle_{t=0}^T$ is locally decomposable
with probability 
$1-o(\vol^{-3\delta/16})$.  By Lemma
\ref{inductionsteplemma} and our inductive assumption that 
(\ref{couplecap}) holds with probability $1-O(\vol^{-\delta/32})$ on 
$\T_k$, 
with probability $1-O(\vol^{-\delta/32})$,
\begin{align*}
\Bigl| \Cap (\T_{k+1})-& \Cap(\T_k)\Bigr| 
\frac{\sqrtvol}{\alpha^{1/2} r}\\  
\leq & \Bigl(\Cap[\LE \langle
X_t \rangle_{t=0}^T] + \close{\T_k}{\LE \langle X_t
\rangle_{t=0}^T}\Bigr) \frac{\sqrtvol}{\alpha^{1/2} r}\\
= & o(\vol^{-\delta/24})
+ \frac{|\Gindex{\ell}| \sqrtvol}{m \alpha^{1/2}r} \\
= & o(\vol^{-\delta/24})  + \frac{|\widetilde \T_{k+1}|
-|\widetilde \T_k|}{\sqrt{m}}
\end{align*}
and hence by induction (\ref{couplecap})
holds with probability $1-O(\vol^{-\delta/32})$.

Finally, we need to show that
the uniform measure $\mu$ and the hitting
measure $\nu$ on $\T_{k+1}$ can be coupled in such a way
that (\ref{uniformhitting}) holds.
Let $\xi$ be chosen according to $\mu$
and $\eta$ according to $\nu$.
Note that 
\begin{equation}
\label{maxdbound}
\dTktwo{\xi}{\eta} \leq \max_{i,j\leq k+1} \dTktwo{x_i}{x_j}.
\end{equation}
We will use (\ref{maxdbound}) as an upper bound whenever
our procedure for constructing the coupling fails.
To make this coupling, 
consider the event $F$ that either $\xi$ and $\eta$
are both in $\T_k$ or 
both are in $\LLE(A_i)$ for some $i\in \Gindex{\ell}$.
If both are in $\T_k$, we wish to use the induction
hypothesis to establish the coupling.  For a random walk $\langle 
Y_u\rangle$ started in the uniform distribution on $G_n$, let $\Told$ 
denote 
the hitting time of $\T_k$ and $\Tnew$ the hitting time of $\LE \langle 
X_t \rangle_{t=0}^T$.
We know that with probability $1-O(\vol^{-\delta/32})$, that 
$X_{\Told}$ 
can be coupled with a uniformly chosen point.  But
$$\P[\Tnew < \Told < \Tnew+r] \leq \frac{\close{\T_k}{\LE \langle 
X_t \rangle_{t=0}^T}}{\Cap[T_{k+1}]}=O(r \vol^{-1/2})
$$
and so
$$
\sum_{x\in \T_k} \bigl| \P[X_{\Told}=x] - \P[X_{\Told}=x \mid \Tnew < 
\Told ]\bigr| \leq O(2^{-r/\mix}) + O(r \vol^{-1/2}).
$$
This means that conditioning on hitting $\T_k$ first has a small enough 
effect that we can still couple $\xi$ and $\eta$ with the claimed 
probability.
 If both $\xi$ and $\eta$
are in $\LLE(A_i)$ for $i\in \good{T}$, then 
$\dTktwo{\xi}{\eta} \leq r=o(\sqrtvol)$. 
We thus need to show that $\P[F^c]=O(\vol^{-\delta/32}).$

To do so, by Lemma \ref{locallydecomplemma},
$\langle X_t \rangle_{t=0}^T$
is locally decomposable with probability $1-o(\vol^{-3\delta/16})$.
On the event $\step{\ell}$,
\begin{align*}
\Cap[\LLE(A_i)] & =\frac{\alpha r^2}{\vol} [1+O(\vol^{-\delta/16})]
\qquad \mbox {\rm and} \\
| \LLE(A_i)| & = \gamma r [1 + O(\vol^{-\delta/12})]
\end{align*}
and so
\begin{equation}
\label{CapdifAi}
\frac{\Cap[\LLE(A_i)]}{|\LLE(A_i)|}=\frac{r \alpha}{\gamma \vol}
\left[1+O\left(\vol^{-\delta/16}\right)\right].
\end{equation}
Likewise, we have already seen that, with probability $1-
O(\vol^{-\delta/32})$, equations
(\ref{couplecap}) and (\ref{couplelength}) hold on $\T_{k+1}$,
so 
\begin{align*}
\Cap(\T_{k+1})& = {|\widetilde \T_{k+1}|}{m} [1+ O(\vol^{-\delta/32})] 
\qquad \mbox{and}\\
|\T_{k+1}|& =|\widetilde \T_{k+1}|(\vol/m)^{1/2} 
[1+O(\vol^{-\delta/32})],
\end{align*}
which gives
\begin{equation}
\label{CapdifTk}
\frac{\Cap[\T_{k+1}]}{|\T_{k+1}|}=\frac{r \alpha}{\gamma \vol}
\left[1+O\left(\vol^{-\delta/32}\right)\right].
\end{equation}
By considering the total variation distance between $\mu$ and
$\nu$, there is a coupling such that
\begin{multline*}
\P[F^c \mid \T_k, \langle X_t \rangle_{t=0}^T]
\leq
\sum_{i\in \Gindex{\ell} \cap \good{T}}
\Bigl| \mu[\LLE(A_i)]-\nu[\LLE(A_i)] \Bigr|  \\ +
\left|\mu(\T_k)-\nu(\T_k)\right| + O\left(\frac{r|\single{T}|}
{|\T_{k+1}|}\right) + O\left(\frac{sT}{r |\T_{k+1}|}\right).
\end{multline*}
To use $\P[F^c]=\E \P[F^c \mid \T_k, 
\langle X_t \rangle_{t=0}^T]$, note that the expected value
of these last two terms is $O(\vol^{-\delta/4})$.
For the other two terms, for any $S$, 
$$|\mu(S)-\nu(S)| \leq \left| \mu(S)-\frac{\Cap(S)}{\Cap(\T_{k+1})}
\right| + \left| \nu(S)- \frac{\Cap(S)}{\Cap(\T_{k+1})}
\right|.
$$
On the event $\step{\ell}$,  
(\ref{CapdifAi}) and (\ref{CapdifTk}) imply that
 for $S=\LLE(A_i)$ or $\T_k$
$$
\left| \mu(S)-\frac{\Cap(S)}{\Cap(\T_{k+1})}
\right|\leq \mu(S) o\left(\vol^{-\delta/32}\right).
$$
Moreover, 
$$
\left| \nu(S)- \frac{\Cap(S)}{\Cap(\T_{k+1})}
\right| \leq \frac{\close{S}{\T_{k+1}\setminus S}}{\Cap[\T_{k+1}]}.
$$
By Lemma \ref{closenesslem},
\begin{equation}
\label{finalexpclosebound}
\frac{
\E [\close{\LLE(A_i)}{\T_{k+1}\setminus 
\LLE(A_i)} 
\mid \T_{k+1}\setminus \LLE(A_i) ]}{\Cap(\T_{k+1})}
\leq 4 \frac{r^2}{\vol}.
\end{equation}
Combining (\ref{finalexpclosebound}) with (\ref{Tkcloseness}), 
we thus obtain
$$
\P[F] \leq \mu(\T_{k+1}) 
o\left(\vol^{-\delta/32}\right)+O\left( \frac{Lr}{\sqrtvol}\right)
+ O\left(\vol^{-\delta/4}\right)
$$
which is $o\left(\vol^{-\delta/32}\right)$ as required.
\end{proof}

\section{Stochastic domination of spanning forests by trees}
\label{treesection}

We now prove Lemma \ref{mainstochasticdomlem}.
The heart of the proof relies on the stochastic domination
of $\F$ by $\extT$.
Suppose that $G$ is a graph with vertex set $V$ and edges $E$.
Let $G_\lambda$ denote the graph that is the extension 
of $G$
formed by adding an additional vertex $\rho$, and from 
every vertex $v\in V$, an edge $(v,\rho)$ of weight $1-\lambda$.
Let $\T_\lambda$ be a weighted spanning tree on $G_\lambda$
as generated by Wilson's algorithm with root vertex $\rho$.  
The graph $\T_\lambda$ induces a forest $\F \subset G$
simply by restricting to edges in $E\cap \T_\lambda$. 
We will show that $\F$ is stochastically dominated by 
the uniform spanning tree $\T$.

An event $A$ is said to be an increasing event on a graph $G$ if
for any subgraph $a\in A$, if $a^\prime$ is 
another subgraph of $G$ formed by adding edges to $a$, then 
$a^\prime\in A$ as well.  We say that an event $A$ is supported on a set 
of edges $E_1$ if determining whether or not $a$ is in $A$ only requires 
looking at the edges $E_1$.  (Equivalently, if $a\in A$, then $a^\prime 
\in A$, where $a^\prime$ is the subgraph whose edges are in both $a$ and 
$E_1$).

\begin{lemma}[Feder and Mihail]
\label{negcorthm}
For increasing events $A$ and $B$ supported on disjoint edge 
sets of $G$,
$
\P[\T \in A \mid \T \in B] \leq \P[\T \in A].
$
\end{lemma}
Let $E_1$ and $E_2$ be disjoint edge sets such that $A$ is supported on 
$E_1$, $B$ is supported on $E_2$, and $E_1\cup E_2=E$.
The case when $|E_2|=1$ was originally proved by Feder and Mihail 
(\cite{FedMih}, Lemma 3.2), and they remark that 
iterating their proof implies that
the general case is also true.  
A proof of the general case appears in the solution to
 Exercise 8.10 in \cite{LyoPerbook}. 

Fix integers $M_{i,j}\geq 1$ for $i,j \leq k$.
Let $B$ denote the 
event that the degree of $\rho$ is at least 2, and let $A$ denote the 
event that there is a path of length exactly $M_{i,j}$ in $G$ from
$x_i$ to $x_j$ for all pairs $i,j \leq k$.
Clearly the event $A$ is increasing and requiring that the path
be in $G$ means that the event
is supported on $E$, the original edge set of $G$.  On the
other hand, the
event $B$ is supported on edges from $G$ to $\rho$.  
As adding edges from $\rho$ to $G$ increases the degree of $\rho$, the 
event $B$ is increasing.  
There are no loops in a forest, 
so the event $A$ implies that
$\dF{x_i}{x_j}=M_{i,j}$ for all $i,j$.

By Lemma 
\ref{negcorthm}, the events $A$ and $B$ are 
negatively correlated.  Moreover, $B^c$ is equivalent to $\F$ 
being a spanning tree.  In particular, conditioned on $B^c$, $\F$ is 
equal to $\T$ in distribution.  This shows that the probability
of having the right lengths in $\F$ is a lower bound for 
the probability of having the right lengths in 
$\T$, i.e., for any collection of values $M_{i,j}<\infty$,
\begin{equation}
\label{stochasticdomeqn}
\P[\dF{x_i}{x_j}=M_{i,j} \, \forall \, i,j \leq k]
\leq 
\P[\dT{x_i}{x_j}=M_{i,j} \, \forall \, i,j \leq k].
\end{equation}  
Using the convention that $\dF{x_i}{x_j}=\infty$
whenever $x_i$ and $x_j$ are in different components of $\F$,
the total variation distance between the joint distribution of
$\dF{x_i}{x_j}$ and the joint distribution of
$\dT{x_i}{x_j}$ is bounded by
$$
\P[\exists \, (i,j): \dF{x_i}{x_j}=\infty].
$$

The key step is now controlling the probability that
$\{x_1, \dots, x_k\}$ are in the same component of $\F$.

\begin{lemma}
\label{intersectlem}
Suppose that $\{ G_n \}$ is a sequence of graphs satisfying the 
assumptions of Theorem \ref{mainthm}.
For any $\eps\in (0,1)$, 
let $\{X_t\}$ and $\{Y_u\}$ be two independent random walks on $G_n$,
with (possibly different) starting points $x$ and $y$, 
and let $T_X$ and $T_Y$ be independent geometric random variables
with mean $\eps^{-2} \sqrtvol$.  Then 
\begin{equation}
\label{intersecteqn}
\P[\LE \langle X_t \rangle_{t=0}^{T_X} \cap \{Y_u\}_{u=0}^{T_Y}
\neq \emptyset] \geq 1-2\eps-a(n,\eps),
\end{equation}
where $a(n,\eps)\rightarrow b(\eps)$ as $n\rightarrow \infty$
and $b(\eps)/\eps \rightarrow 0$ as $\eps \rightarrow 0$.
\end{lemma}

\begin{proof} 
Let $s$ and $r$ be as in Table \ref{deftable}.
For any set $S\subset G_n$, 
let $T_S$ denote the time $\{Y_u\}_{u=0}^\infty$ 
first hits $S$.  If $\Cap(S) \geq r \eps/\sqrtvol$,
then
\begin{align*}
\P[T_Y < T_S] & \leq \sum_{i=0}^\infty
\P[ir < T_Y \leq  (i+1)r \mid T_S, T_Y > ir ] \P[T_S, T_Y >ir]\\
& \leq \sum_{i=0}^\infty \frac{\eps^2 r}{\sqrtvol}
\P[T_S > ir] \P[T_Y > ir]\\
& \leq \frac{\eps^2 r}{\sqrtvol} \sum_{i=0}^\infty
\left( 1-\frac{r \eps}{2\sqrtvol}\right)^i \\
& = \frac{\eps^2 r}{\sqrtvol} \frac{2 \sqrtvol}{r \eps}\\
& = 2 \eps.
\end{align*}
The proof thus reduces to finding $a(n,\eps)$ such that
$$
\P\left[
\Cap [\LE \langle X_t \rangle_{t=0}^\infty] > \frac{r 
\eps}{\sqrtvol}\right]
\geq 1-a(n,\eps).
$$

Let $T=(2\eps/\alpha) \sqrtvol$, where $\alpha$ is as in
(\ref{alphagammadef}).
Since $T_X$ is geometric with mean 
$\eps^{-2} \sqrtvol$, we have
$
\P[T_X < T ] \leq 2 \eps^{3}/\alpha.
$
Restricting to the event $\{T_X \geq T\}$, 
we will now view the path 
$\langle X_t \rangle_{t=0}^{T_X}$ as being in two
parts: an initial run of length $T_X-T$ and a final tail
of length $T$.  The idea of the proof is that if the initial segment
has a large capacity, then because the tail is short, it probably
misses enough of the initial segment that the capacity remains high.
Conversely, if the initial segment has low capacity,
then the tail will probably survive and is long enough to be of high
capacity itself.  

More formally, consider two possibilities 
based on whether or not the event
\begin{equation}
\label{initialcapeqn}
\{\Cap [\LE \langle X_t 
\rangle_{t=0}^{T_X-T}] > 2\eps r/\sqrtvol\}
\end{equation}
occurs.
Let $\langle \eta(k)\rangle$
be an increasing sequence such that
$\LE \langle X_t \rangle_{t=0}^{T_X-T}=
\langle X_{\eta(k)} \rangle$.
When (\ref{initialcapeqn}) holds, 
let $M$ be the smallest number such that 
$$\Cap[\langle X_{\eta(k)} \rangle_{k\leq M}] 
> \frac{\eps r}{\sqrtvol}.$$
Denote this initial segment by 
$U:= \langle X_{\eta(k)} \rangle_{k\leq M}.$
Because capacity is subadditive, $M\geq \eps \sqrtvol$ and
$$
\frac{\eps r}{\sqrtvol} < \Cap U < \frac{\eps r}{\sqrtvol}+ 
\frac{r}{\vol}<\frac{2\eps r}{\sqrtvol}.
$$
Subadditivity of capacity also implies that 
$T_X-T-M> \eps \sqrtvol-1$, which for large enough $n$ 
is greater than $s$.
But
\begin{align*}
\P[\{X_t\}_{t=T_X-T}^{T_X-T+s} \cap U \neq \emptyset]
& \leq \E \Bigl|\{X_t\}_{t=T_X-T}^{T_X-T+s} \cap 
\{X_t\}_{t=0}^{T_X-T-\eps \sqrtvol} \Bigr| \\
&\leq \frac{2s \E T_X}{\vol}=\frac{2s}{\eps^2 \sqrtvol}.
\end{align*}  
Subdividing the final $T-s$ steps of the walk into pieces of length
$r$, using the fact that $X_{T_X-T+s}$ is close to uniform, and
considering the expected number of those pieces that intersect $U$
gives
$$
\P[\{X_t\}_{t=T_X-T+s}^{T_X} \cap U \neq \emptyset] \leq \frac{2 \eps 
r}{\sqrtvol}
\frac{2T}{r}
 = \frac{8 \eps^2}{\alpha}.
$$
If $\{X_t\}_{t=T_X-T}^{T_X} \cap U= \emptyset$, then $U$ survives
loop-erasure, so this shows that, 
conditioned on (\ref{initialcapeqn}),
with probability $1-a(n,\eps)$ for $a(n,\eps)$ of the suitable form, 
$U$ survives loop-erasure and yields the 
desired capacity for $\LE\langle X_t\rangle$.

When (\ref{initialcapeqn}) fails, we have
$\Cap[\LE \langle X_t \rangle_{t<T_X-T}]\leq 2\eps r/\sqrtvol$.  The 
probability
of a segment of length $T$ started at uniform intersecting
this initial piece is then $O(\eps^{2})$ by the same argument as
before.  Adding a buffer of $s$ steps to get close to a uniform position,
we thus have
$$\P[\{X_t\}_{t >T_X-T+s} \cap \LE \langle X_t \rangle_{t=0}^{T_X-T} = 
\emptyset]=1-O(\eps^{2})-(1/2)^{\lfloor s/\mix \rfloor}.$$
Moreover, $T$ is small 
enough such that the expected number of loops longer than $\mix$ within 
the final $T-\mix$ steps is bounded by $4 \eps^2/\alpha^2$.  Thus with 
probability $1-f(n,\eps)$, where $f(n,\eps)=O(\eps^2)+(1/2)^{\lfloor
s/\mix \rfloor}$,
$\LE \langle X_t \rangle_{t=0}^{T_X}$ contains
$W=\LLE \langle X_t \rangle_{t=T_X-T+s}^{T_X}$,
and so
\begin{align*}
\P\left[ \Cap(\LE \langle X_t \rangle_{t=0}^{T_X})> \frac{r \eps}
{\sqrtvol}\right]
& \geq \P\left[ \Cap(W)
> \frac{r \eps} {\sqrtvol}\right]- f(n,\eps)\\
& = \P\left[ \Cap(W)
> \frac{\alpha}{2}\frac{T r} {\vol}\right]- f(n,\eps).
\end{align*}
By again breaking the final $T-s$ steps up into runs of length
$r$, using the concentration of capacity about its mean on 
such a run, and bounding the closeness between these runs
as in the proof of Lemma \ref{couplinglemma}, this final
probability is of the form $a(n,\eps)$ as required, with
$b(\eps)=O(\eps^2)$.
 \end{proof}

\begin{proof}[Proof of Lemma \ref{mainstochasticdomlem}]
As mentioned above, we need to bound
\begin{equation}
\label{TVneedeqn}
\P[\exists \, (i,j): \dF{x_i}{x_j}=\infty],
\end{equation}
and in particular show that it is less than $\eps$
for sufficiently large $n$ and $L$.
But using Lemma \ref{intersectlem} gives
\begin{align*}
\P[\exists \, (i,j): \dF{x_i}{x_j}=\infty] & \leq 
\frac{k^2}{2} \P[\dF{x_1}{x_2}=\infty] \\
& \leq \frac{k^2}{2} \frac{2}{L^{1/2}} + a(n, L^{-1/2}),
\end{align*}
where $a(n,x)$ is as in Lemma \ref{intersectlem}.
Taking $L=k^4 \eps^{-2}$ bounds
(\ref{TVneedeqn}) by $\eps + a(n,\eps k^{-2})$,
which for small enough $\eps$ and large enough $n$
yields the required bound.
\end{proof}

\section{Constants on the torus}
\label{torussection}

This section is devoted to proving Theorem \ref{torusthm}.  
To see that $\Z_n^d$ satisfies the hypotheses of Theorem \ref{mainthm}, 
note that for simple random walk on $\Z_n^d$ with holding probability 
$1/2$, $|\P_o[X_t=o]-n^{-d}|\leq C 
t^{-d/2}$ for a suitable constant $C$ (see e.g., \cite{AldFil}, Chapter 
5). 
This means that the
only way in which Theorem \ref{torusthm} is not a special case of Theorem 
\ref{mainthm} is that there is a single rescaling 
constant $\beta$ rather than a sequence of constants $\beta_n$ that 
(possibly) depend on $n$.  Thus, what we need to show in this 
section is that $\lim \beta_n$ exists.  We will 
do so by giving an expression for the limit.

\begin{lemma}
Let $\{\hatY_t\}$, $\{\hatZ_u\}$, and $\{\hatW_v\}$ be independent simple 
random 
walks 
on $\Z^d$, all starting at the origin.  Let $G_n=\Z_n^d$, and take 
$\alpha=\alpha(n)$ and $\gamma=\gamma(n)$ to be as in 
(\ref{alphagammadef}). Then for $d \geq 5$,
\begin{equation}
\label{torusgamma}
\limn \gamma(n) = \P[ \LE \langle \hatY_t\rangle_{t=0}^\infty \cap 
\{\hatZ_u\}_{u=1}^\infty = \emptyset] 
\end{equation}
\begin{multline}
\label{torusalpha}
\limn \alpha(n)= \\
\P\left[ \LE \langle \hatY_t\rangle_{0}^\infty \cap
\{\hatZ_u\}_{1}^\infty = \emptyset, \left(\LE \langle 
\hatY_t\rangle_{0}^\infty \cup \LE \langle 
\hatZ_u\rangle_{1}^\infty\right) 
\cap \{\hatW_v\}_{1}^\infty = \emptyset\right].
\end{multline}
\end{lemma}

\begin{proof} 
Take $s, q,$ and $r$ as in Table \ref{deftable}.  To
understand scales, recall that on $\Z_n^d$, $\mix$ is on the order of
$n^2$, meaning that $s,q$ and $r$ are on the order of $n^{(d+12)/8}$,
$n^{(d+4)/4}$, and $n^{(3d+4)/8}$ respectively.

Let $U\subset [2s+1,r-s]$ be the set of times $t$ in the interval 
$[2s+1,r-s]$ 
that are locally retained (see Definition \ref{localretaineddef}).
As the probability of each time being locally retained is constant
on $[2s+1,r-s]$,
\begin{align*}
\gamma r & =\E | \LLE \langle X_t \rangle_{t=2s+1}^{r-s}| \\
& = \sum_{t=2s+1}^{r-s} \P[t\in U]\\
& = (r-3s) \P[q \in U].
\end{align*}
As $s=o(r)$, proving (\ref{torusgamma}) reduces to 
computing $\lim_{n\rightarrow \infty} \P[q \in U]$.
On the torus in dimension $d\geq 5$, the expected number
of loops of length greater than $n^{7/4}$ in a run of length
$r$ is bounded by 
\begin{align*}
\sum_{i=0}^r \sum_{j=i+n^{7/4}}^r \P[X_i=X_j] & 
\leq \sum_{k=n^{7/4}}^r r \P_o[X_k=o]\\
& = O\left( r^2 n^{-d} + r \left(n^{7/4}\right)^{1-d/2}\right) \\
& = O\left( n^{(4-d)/4} + n^{(18-4d)/8} \right).
\end{align*}
For $d\geq 5$, this expression is $o(1)$.
Moreover, Lemma \ref{cutpointlem1} and the proof of Corollary
\ref{cutptcor} imply that
$$\P[\exists\, T\in [n^{7/4}, n^{9/5}] : 
\langle X_t \rangle_{t=T-n^{7/4}}^{T} \cap
\langle X_t \rangle_{t=T+1}^{T+n^{7/4}} = \emptyset]
=1-o(1).$$
Combining these two facts, the probability that 
$\langle X_t \rangle_{t=0}^s$ has a cutpoint
in the time interval $n^{7/4}\leq t \leq n^{9/5}$ is
$1-o(1)$.
The importance of having cutpoints
is that whether or not the point $X_j$ survives loop-erasure can
be determined from only considering what happens
between two cutpoints, one at a time before $j$, and one at a
time after $j$.
As $9/5<2$, with probability $1-o(1)$,
any run of length $n^{9/5}$ 
remains inside a cube of edge length $n$, and in particular
does not see the difference between the torus and the full
lattice $\Z^d$.  Combining these facts,
\begin{align}
\P[q \in U] & = \P [\LE \langle X_t \rangle_{t=-s}^0 \cap 
\{X_t\}_{t=1}^s = \emptyset]  \notag \\
& = \P [ \LE \langle X_t \rangle_{t=-n^{9/5}}^0 \cap
\{X_t\}_{t=1}^{n^{9/5}} = \emptyset] +o(1) \notag \\
\label{fulllatticeeqn}
& = \P [ \LE\langle \hatY_t \rangle_{t=0}^{n^{9/5}} 
\cap \{\hatZ_u\}_{u=1}^{n^{9/5}} = \emptyset] +o(1).
\end{align}
Note that
if $\left(\LE \langle \hatY_t \rangle_{t=0}^\infty \cap 
\{\hatZ_u\}_{u=1}^\infty\right) \neq
\left(\LE \langle \hatY_t \rangle_{t=0}^{n^{9/5}} \cap
\{\hatZ_u\}_{u=1}^{n^{9/5}}\right)$, then
either $\{\hatY_t\}_{t>n^{9/5}}\cap \{\hatZ_u\}_{u>0}\neq \emptyset$
or $\{\hatY_t\}_{t>0}\cap \{\hatZ_u\}_{u>n^{9/5}}\neq \emptyset$.
But a first moment argument shows that
$$
\P\left[\left(\{\hatZ_u\}_{u=n^{9/5}}^\infty \cap 
\{\hatY_t\}_{t=1}^\infty\right)
=\left(\{\hatZ_u\}_{u=1}^\infty \cap \{\hatY_t\}_{t=n^{9/5}}^\infty\right)
=\emptyset\right]=1-o(1),
$$
which completes the proof of (\ref{torusgamma}).

To prove (\ref{torusalpha}), let 
$S=\LLE \langle X_t \rangle_{t=2s+1}^{r-s}$. 
As $\alpha= r^{-2} \E \Cap S \vol$, we need to compute $\E \Cap S$.
 For a simple random walk 
$\{Y_k\}_{k=0}^\infty$ on $\Z_n^d$, let $\tau_S=\inf\{k\geq 0: Y_k \in 
S\}$, and 
$T_S=\inf\{k \geq 1: 
Y_k \in S\}$.   Considering the time reversal and
again letting $U\subset [2s+1,r-s]$ denote the locally 
retained times,
\begin{align*}
\Cap S = & \sum_{j\in U} \sum_{k=0}^r \P_{\pi} [\tau_S=k, X_j=Y_k] \\
 = & \sum_{j\in U} \sum_{k=0}^r \sum_{z\in \Z_n^d} \P_z[\tau_S=k,
 X_j=Y_k] n^{-d}\\
 = &\sum_{k=0}^r \sum_{j \in U} \sum_{z \in {\Z_n^d}} \P_{X_j} [T_S>k, 
Y_k=z] n^{-d} \\
= & \sum_{k=0}^r \sum_{j \in U} n^{-d} \P_{X_j} [T_S>k].
\end{align*}

Let $\1_U(\cdot)$ be an indicator function for $U$ 
and let $W$ denote the event $\{\LE \langle X_t \rangle_{t=-s}^0 \cap 
\{X_t\}_{t=1}^{s} = \emptyset\}$.  Then
\begin{align*}
\E \Cap S = & \sum_{k=0}^r \sum_{j=2s+1}^{r-s}
\E \left[\E \left[ \1_U(j) \P_{X_j}(T_S>k) 
\mid \langle X_t \rangle_{t=0}^r 
\right]\right] n^{-d} \\
= & \sum_{k=0}^r \sum_{j=2s+1}^{r-s} n^{-d} 
\E\left[ \P [ W,
\LLE \langle X_t \rangle_{2s+1-j}^{r-s-j} \cap \{Y_u\}_{1}^k = 
\emptyset \mid Y_0=X_0]\right].
\end{align*}
There are fewer than $3rs$ terms in which
$k < s$ or $j \notin [3s, r-2s]$,
each of which is bounded by $n^{-d}$.  The sum of these terms thus
contributes at most $3rsn^{-d}$, which is of a lower order than 
$\E\Cap S$ (which is on the order of $r^2 n^{-d}$). 
As $\E \Cap S= \alpha r^2 n^{-d}$, 
it thus suffices to show that
for the $r^2(1+o(1))$ terms
with $k\geq s$ and $j \in [3s, r-2s]$, we uniformly obtain
\begin{align*}
\limn 
& \E \left[\P [W, \LLE \langle X_t \rangle_{2s+1-j}^{r-s-j} 
\cap \{Y_u\}_{1}^k = \emptyset \mid Y_0=X_0]\right] \\
= \P & \left[ \LE \langle Y_t\rangle_{0}^\infty \cap
\{Z_u\}_{1}^\infty = \emptyset, \left(\LE \langle 
\hatY_t\rangle_{0}^\infty \cup \LE \langle 
\hatZ_u\rangle_{1}^\infty\right) 
\cap \{\hatW_v\}_{1}^\infty = \emptyset\right].
\end{align*}
As before, since we are only running the walk for times on the order 
of $r$, the probability that there are no loops of length 
longer than $n^{7/4}$ is $1-o(1)$.  We then convert from a statement on 
the torus to one on the full lattice exactly as before.
The convergence is uniform because 
the analogous conversions to (\ref{fulllatticeeqn})
rely only on the existence of these cutpoints.
\end{proof}

\section{Expanders, Hypercubes, and Proof of Theorem \ref{generalCRTthm}}
\label{gensection}

We stated in Section \ref{introduction} that sequences of expander graphs 
satisfy the assumptions of Theorem \ref{mainthm}.  This is immediate from 
the fact that, for a sequence of expanders, there exist $C>0$ and 
$\lambda<1$ such that the bound 
$|\P_o[X_t=o]-\vol^{-1}| \leq C \lambda^t$ holds for the entire sequence.  

We likewise claimed that (\ref{superstrongtran}) applies on the 
hypercubes $\Z_2^n$.  To see this, consider two 
random walks $\{X_t\}$ and $\{Y_u\}$ with different starting points, both 
run in continuous time 
with rate 1.  The expected amount of time of intersection is the expected 
number of intersections for two discrete time walks.
In continuous time, 
\begin{equation}
\label{hypercubeeqn}
\P[X_t=Y_u]\leq \left(\frac{1}{2}\right)^n \left(1+e^{-2(t+u)/n} 
\right)^{n-1} \left(1-e^{-2(t+u)/n}\right).
\end{equation}
For $t+u \leq n^{1/4}$,
bound (\ref{hypercubeeqn}) by 
 $(1 - \exp[-2(t+u)/n])/2<2(t+u)/n$.
As $e^{-x} \leq 1-x/2$ for $0\leq x\leq 1$, we obtain the bound
$$
\frac{1}{2} \left(1+ \exp\left[-\frac{2(t+u)}{n}\right]\right)
< 1-\frac{t+u}{n} < \exp\left[-\frac{t+u}{n}\right]
$$
for $n^{1/4}< t+u \leq n/2$,
and then finally
$$\left(\frac{1}{2}\right)^n \left(1+e^{-2(t+u)/n}
\right)^{n-1} \leq \left(\frac{1}{2}\right)^n\left(
1+e^{-1}\right)^n$$ 
for $(t+u)>n/2$ yields
$$
\int_0^r \int_0^r \P[X_t=Y_u] \, dt \, du =o(1) \, ,
$$
which in turn says that the expected time of overlap of the two paths is 
$o(1)$.  As we ran the continuous time walks for time $r \gg q$, the 
expected number of intersections in the first $q$ steps of the discrete
walks is also $o(1)$.

Turning now to the proof of Theorem \ref{generalCRTthm}, note that 
Theorem \ref{generalCRTthm} differs from Theorem \ref{mainthm}
in two ways: first, we need to show that assumption 
(\ref{superstrongtran}) allows us to omit
the hypothesis that  $\{x_1, \dots, x_k\}
=\{x_1^{(n)}, \dots, x_k^{(n)}\}$ 
are chosen uniformly, 
and second, we need to show that $\lim \beta_n=1$.

To show that we can choose $\{x_1, \dots, x_k\}$
freely under assumption (\ref{superstrongtran}),
let $\T_0=\{\rho\}$, and for $1\leq i \leq k$,
let $\{X_t^i\}_{t=0}^\infty$ be i.i.d., weighted
random walks on $\extG{n}{L}$
with $X_0^i=x_i$.  Let $T_i=\min\{t\geq 0:
X_t^i\in \T_{i-1}\}$, and take $\T_i=\T_{i-1}\cup \LE
\langle X_t^i \rangle_{t=0}^{T_i}$.

Let $s$ be as in Table \ref{deftable}.
As the $L^\infty$ distance between the distribution of $X_t^i$
and the uniform is $o(\vol^{-1})$ for $t\geq s$,
we only need to show that $\P[T_k <s ]=o(1)$.
The main concern is intersections that might occur
from $x_k$ being close to $\{x_1, \dots, x_{k-1}\}$.
After running for $s$ steps, the $X_t^i$ are close
to uniform.   Using the expected number
of intersections between $\{X_t^k\}_{t<s}$
and $\{X_t^i\}_{t\geq s}$ to bound the probability
of such an intersection occurring, we obtain
\begin{equation}
\label{endeqn}
\P[T_k < s]
\leq k \max_{i<k} \P\left[\{X_t^i\}_{t<s}
\cap \{X_t^k\}_{t<s}  \neq \emptyset \right] + \frac{2 s \E 
\sum_{i=1}^{k-1} T_i}{\vol}.
\end{equation}
But assumption (\ref{superstrongtran}) is the fact that  
$\P\left[\{X_t^i\}_{t<s} \cap \{X_t^k\}_{t<s}  \neq \emptyset 
\right]=o(1)$, so as $\E T_i \leq L \sqrtvol$, equation
(\ref{endeqn}) is exactly what we need.

For the second part, assumption (\ref{superstrongtran})
implies that any point is a local cutpoint with probability $1-o(1)$, so 
$\gamma=1-o(1)$.  Likewise, (\ref{superstrongtran}) implies that the 
probability of a run of length $r$ 
intersecting $\LLE(A_i)$ more than once, even conditioned on there being 
an intersection, is $o(1)$.  This means that the probability of an 
intersection is, up to a factor of $1+o(1)$, the same as the expected 
number of intersections, and so we also 
have $\alpha=1-o(1)$.  In particular, $\beta_n=1+o(1)$, which completes 
the 
proof of Theorem \ref{generalCRTthm}.

\section{Further questions}
\label{questionsection}

Although the main results of this paper give a good
picture of the scaling limit of UST on many graphs,
there are still a number of questions that remain.

\begin{enumerate}
\item
Is the UST on the complete graph in some sense smaller than on any other 
vertex transitive graph?  More precisely, if $\{G_n\}$ are vertex 
transitive, and $x$ and $y$ are uniformly chosen from $G_n$, is there a 
constant $C$ such that
$$
\P[\dT{x}{y} > \lambda \sqrtvol ] \geq 
\exp\left[-C \frac{\lambda^2}{2}\right] (1+o(1))?
$$
Benjamini and Kozma \cite{BenKoz03} asked an averaged form of this 
question, 
asking if 
$\E \dT{x}{y} \geq C \sqrtvol$ holds.
 \item
Theorems \ref{torusthm}-\ref{generalCRTthm} only prove
that the scaling limit of the UST is the Brownian
CRT in the sense that the finite dimensional
distributions converge.  
Does this convergence also hold in a stronger topology?
\item
Our theorems do not apply to the
torus $\Z_n^4$ because
$\mix$ and $\sqrtvol$ are on the same order of magnitude
in dimension 4.
After taking into account a logarithmic correction factor,
the scaling limit of LERW on $\Z^4$, however, is
still Brownian motion \cite{Lawbook}.  As discussed in
\cite{BenKoz03}, heuristics suggest that 
$\E \dT{x}{y}$ is on the order of $n^2 \log^{1/6}n$.  
If so,
what is the limiting distribution of $\dT{x}{y}$?
What is the scaling limit of the UST on $\Z_n^4$?
\item
In this paper, we have focused on the intrinsic geometry
of the UST, discussing distances in the UST.  
We can also ask about the existence of a scaling limit in the
extrinsic geometry induced by embedding our graph in the torus of side 
length 1.  The path of the LERW is asymptotically dense in this embedding,
but lifting to the universal cover $\R^d$ and dividing lengths by 
$n^{d/4}$ (the square root of the typical length of a path), 
Theorem \ref{torusthm} suggests the following scaling limit for the 
lifted UST:  first, the LERW from $x$ 
to $y$ lifts to a Brownian motion on $\R^d$ run for a random amount of 
time $T$ that is Rayleigh distributed.  Lifting the partial spanning tree 
defined by $k$-points in this way, 
 we obtain an embedding in $\R^d$ of the first $k$ steps of the Poisson 
line breaking constructing for the Brownian CRT, where an edge length of 
length $\ell$ in the CRT corresponds to a Brownian path run for time 
$\ell$. (This is a version of Le-Gall's Brownian snake.) 
However, establishing this picture requires further work.
\end{enumerate}

\noindent{\bf Acknowledgement}
We  are grateful to Jim Pitman for sharing his conjecture
relating the scaling limit of UST on the torus to the Brownian
Continuum random tree, and for corrections to an earlier version of this 
manuscript.
We also thank Ron Peled, G{\a'a}bor Pete, Jason Schweinsberg, and
B{\'a}lint Vir{\'a}g for useful comments and discussions.

\bibliographystyle{amsplain}
\bibliography{../refs} 

\end{document}